\documentclass[10pt]{amsart}
\usepackage{amssymb}
\usepackage{dsfont}
\usepackage{amscd}
\usepackage[mathscr]{euscript} 
\usepackage[all]{xy}
\usepackage{url}
\usepackage{stmaryrd}
\usepackage{comment}
\usepackage[retainorgcmds]{IEEEtrantools}
\usepackage{hyperref}
\usepackage{graphicx}
\usepackage{mathtools}
\usepackage[dvipsnames]{xcolor}
\usepackage{centernot}
\usepackage{marvosym}
\usepackage{tikz}
\usepackage{tikz-cd}
\usepackage{xspace}

\usepackage{tkz-fct}
\usepackage{soul}
\usepackage{lineno}
\usepackage[T1]{fontenc}
\usepackage[utf8]{inputenc}

\usepackage[english]{babel}
\usepackage[autostyle]{csquotes}
\usepackage[shortlabels]{enumitem}
\usepackage[backend=biber, style=alphabetic]{biblatex}
\DeclareMathOperator{\Th}{Th}

\title[Ranks in Ellis semigroups and model theory] {Ranks in Ellis semigroups and model theory}


\author[A. Codenotti]{Alessandro Codenotti$^{\ast}$}
\thanks{$^{\ast}$ The first author is funded by the Deutsche Forschungsgemeinschaft (DFG, German Research Foundation) under Germany’s Excellence Strategy EXC 2044–390685587, Mathematics Münster: Dynamics– Geometry–Structure and by CRC 1442 Geometry: Deformations and Rigidity.}
\address{$^{\ast}$ Institut f\"{u}r Mathematische Logik und Grundlagenforschung, Universit\"{a}t M\"{u}nster, Einsteinstrasse 62, 48149 M\"{u}nster, Germany}
\email{acodenot@uni-muenster.de}

\author[D. M. HOFFMANN]{Daniel Max Hoffmann$^{\dagger}$}
\thanks{$^{\dagger}$SDG. The second author is supported by the European Union’s Horizon research and innovation
programme under the MSCA project no. 101063183,
by the National Science Centre (Narodowe Centrum Nauki, Poland) 
grant no. 2021/43/B/ST1/00405, and by the UW IDUB PSP no. 501-D110-20-3004310.}
\address{$^{\dagger}$
Institut f\"{u}r Geometrie\\
Technische Universit\"{a}t Dresden\\
Dresden\\
Germany
\newline {\em and}
\newline
Instytut Matematyki\\
Uniwersytet Warszawski\\
Warszawa\\
Poland}
\email{daniel.max.hoffmann@gmail.com}
\urladdr{{https://sites.google.com/site/danielmaxhoffmann/}}

\date{\today}

 \DeclareMathOperator{\aut}{Aut} \DeclareMathOperator{\id}{id}

\DeclareMathOperator{\tp}{tp}

\DeclareMathOperator{\ddf}{DF}\DeclareMathOperator{\dcf}{DCF}\DeclareMathOperator{\scf}{SCF}

\DeclareMathOperator{\CB}{CB}

\DeclareMathOperator{\res}{res}

\DeclareMathOperator{\osc}{osc}
\DeclareMathOperator{\Ord}{Ord}
\DeclareMathOperator{\EL}{EL}
\DeclareMathOperator{\bR}{\beta R}

\DeclareMathOperator{\Nn}{\mathbb{N}}

\DeclareMathOperator{\Cc}{\mathbb{C}}

\newtheorem{theorem}{Theorem}[section]
\newtheorem{prop}[theorem]{Proposition}
\newtheorem{lemma}[theorem]{Lemma}
\newtheorem{cor}[theorem]{Corollary}

\theoremstyle{definition}
\newtheorem{definition}[theorem]{Definition}
\newtheorem{example}[theorem]{Example}
\newtheorem{remark}[theorem]{Remark}
\newtheorem{question}[theorem]{Question}

\theoremstyle{remark}
\newtheorem*{theorem*}{Theorem}
\newtheorem*{cor*}{Corollary}

\theoremstyle{definition}
\theoremstyle{definition}

\theoremstyle{definition}

\theoremstyle{remark}

\makeatletter
\providecommand*{\cupdot}{%
  \mathbin{%
    \mathpalette\@cupdot{}%
  }%
}
\newcommand*{\@cupdot}[2]{%
  \ooalign{%
    $\m@th#1\cup$\cr
    \sbox0{$#1\cup$}%
    \dimen@=\ht0 %
    \sbox0{$\m@th#1\cdot$}%
    \advance\dimen@ by -\ht0 %
    \dimen@=.5\dimen@
    \hidewidth\raise\dimen@\box0\hidewidth
  }%
}

\providecommand*{\bigcupdot}{%
  \mathop{%
    \vphantom{\bigcup}%
    \mathpalette\@bigcupdot{}%
  }%
}
\newcommand*{\@bigcupdot}[2]{%
  \ooalign{%
    $\m@th#1\bigcup$\cr
    \sbox0{$#1\bigcup$}%
    \dimen@=\ht0 %
    \advance\dimen@ by -\dp0 %
    \sbox0{\scalebox{2}{$\m@th#1\cdot$}}%
    \advance\dimen@ by -\ht0 %
    \dimen@=.5\dimen@
    \hidewidth\raise\dimen@\box0\hidewidth
  }%
}
\makeatother

\def\Ind#1#2{#1\setbox0=\hbox{$#1x$}\kern\wd0\hbox to 0pt{\hss$#1\mid$\hss}
\lower.9\ht0\hbox to 0pt{\hss$#1\smile$\hss}\kern\wd0}

\def\notind#1#2{#1\setbox0=\hbox{$#1x$}\kern\wd0
\hbox to 0pt{\mathchardef\nn=12854\hss$#1\nn$\kern1.4\wd0\hss}
\hbox to 0pt{\hss$#1\mid$\hss}\lower.9\ht0 \hbox to 0pt{\hss$#1\smile$\hss}\kern\wd0}

\begin{document}

\newcommand{\ov}{\overline}
\newcommand{\FC}{\mathfrak{C}}

\newcommand{\twoc}[3]{ {#1} \choose {{#2}|{#3}}}
\newcommand{\thrc}[4]{ {#1} \choose {{#2}|{#3}|{#4}}}
\newcommand{\Kk}{{\mathds{K}}}

\newcommand{\dlog}{\mathrm{ld}}
\newcommand{\ga}{\mathbb{G}_{\rm{a}}}
\newcommand{\gm}{\mathbb{G}_{\rm{m}}}
\newcommand{\gaf}{\widehat{\mathbb{G}}_{\rm{a}}}
\newcommand{\gmf}{\widehat{\mathbb{G}}_{\rm{m}}}
\newcommand{\gdf}{\mathfrak{g}-\ddf}
\newcommand{\gdcf}{\mathfrak{g}-\dcf}
\newcommand{\fdf}{F-\ddf}
\newcommand{\fdcf}{F-\dcf}
\newcommand{\mw}{\scf_{\text{MW},e}}

\newcommand{\BC}{{\mathbb C}}

\newcommand{\CC}{{\mathcal C}}
\newcommand{\CG}{{\mathcal G}}
\newcommand{\CK}{{\mathcal K}}
\newcommand{\CL}{{\mathcal L}}
\newcommand{\CN}{{\mathcal N}}
\newcommand{\CS}{{\mathcal S}}
\newcommand{\CU}{{\mathcal U}}
\newcommand{\CF}{{\mathcal F}}
\newcommand{\CP}{{\mathcal P}}
\newcommand{\CI}{{\mathcal I}}
\newcommand{\SL}{{\mathrm{SL}}}

\begin{abstract}
    We slightly generalize a notion of rank introduced by Glasner and Megrelishvili, which captures the oscillations of elements of Ellis semigroups, so that it can be applied to any compact Hausdorff space instead of being limited to the metric case. Then, we relate this rank to classical dividing lines in the model-theoretic stability hierarchy. For example, that the rank is ordinal-valued if and only if the background theory is NIP.
\end{abstract}

\maketitle

\section{Introduction}
In January 2023, the first author was giving a talk about (metrizable) tame dynamical systems at the \textit{Bonn-Münster-Düsseldorf GeSAMT - Gemeinsames Seminar Algebra und Modelltheorie}. After providing the definition of the rank from  \cite{GM} (cf. \cite{KechrisLouveau}, \cite{HOR}), the second author, being in the audience, asked how is this rank related to model theory, having in mind that the Cantor-Bendixson rank is the Morley rank, and noticing some similarities of the rank from \cite{GM} to the Cantor-Bendixson rank. This manuscript is a common effort by both authors to answer the question.

For us, a compact space is not necessarily Hausdorff. On the other hand, we always assume a dynamical system to have compact Hausdorff space as one of its main ingredients:

\begin{definition}
    A pair $(X,G)$, where $X$ is a topological space and $G$ is a topological group, is called \emph{dynamical system} if $X$ is compact Hausdorff and $G$ acts continuously via homeomorphisms on $X$.
\end{definition}

Let $(X,G)$ be a dynamical system (we keep this assumption until the end of this paper).
The typical example we had in mind is the space of types $S(\FC)$ of a monster model $\FC$ (i.e. $\kappa$-saturated and strongly $\kappa$-homogeneous model, for some big $\kappa$, e.g. exists for every uncountable regular $\kappa$ such that $\kappa=2^{<\kappa}$) of an $\CL$-theory $T$. 
The first natural group action is the action of $\aut(\FC)$ on $S(\FC)$
(in model theory actions on spaces of types related to definable groups are studied as well \cite{ArtemPierre}, other possibilities include to study spaces of specific families of types or even Keisler measures).
As we know, the topological Cantor-Bendixson rank on the space of types $S(\FC)$ is the model-theoretic Morley rank.
Here we ask, what is the $\beta$-rank introduced in \cite{GM} when we consider the dynamical system $(S(\FC),\aut(\FC))$. However the original $\beta$-rank from \cite{GM} was defined for metric spaces, while we want to study dynamical systems coming from model theory and the topological spaces there are not necessarily metric. We begin Section \ref{section: definitions and first observations} by slightly generalizing the definition of $\beta$-rank to work for arbitrary compact Hausdorff spaces.  After the first definitions and observations, we begin to investigate the relationship between $\beta$-rank and model theoretic dividing lines. In particular we show that a theory $T$ is stable if and only if $\beta(S_x(\FC),\aut(\FC))=0$ (cf. Definition \ref{definition: beta rank}), where $\FC\models T$ is a monster model, in Proposition \ref{prop: rank 0 iff stable}.
In other words, a theory $T$ is stable if and only if every element of the Ellis semigroup $E(S_x(\FC),\aut(\FC))$ is a continuous map. The last thing was already known, but perhaps never stated without additional technical assumptions like $\aleph_0$-categoricity (see \cite[Section 5]{WAPstable}).

In Section \ref{section: o-minimal}, we study o-minimal theories and show that $\beta(S_x(\FC),\aut(\FC))=1$ for $T$ being o-minimal, $\FC\models T$ being a monster model and $x$ being a single variable. This is new, but based on techniques from \cite{codenotti2023}. We also provide an example which is not o-minimal but with $\beta$-rank equal to $1$,
thus showing that the other implication does not hold.

Section \ref{section: NIP theories} is the main part of this paper, and is devoted to relating with each other three properties: tameness of $(S_x(M),\aut(M))$, the property that $\beta(S_x(M),\aut(M))$ has an ordinal value and ``NIP"-ness of a theory $T$ of which $M$ is a model.
In Corollary \ref{cor: tame vs defined rank}, we show that $(X,G)$ is tame if and only if $\beta(X,G)<\infty$.
This is a general result from topological dynamics which was already known in the metric case. 
The definition of tameness of a dynamical system originated in the work of K\"ohler, where tame dynamical systems are called \emph{regular} \cite{Kohler}. Tameness is now considered an important property and has been studied by various authors, see: \cite{KerrLi}, \cite{Glasner1}, \cite{Glasner2},\cite{Glasner3}, 
\cite{Huang}, \cite{GMtrees} and \cite{GM23}.
It is known that, for a metric dynamical system $(X,G)$, tameness is equivalent to the statement that every function $f$ in the Ellis semigroup $E(X,G)$ is Baire class 1. We could not follow this path, as we are not assuming metrizability and hence 
for us, Definition 3.1 from \cite{GM23}, in which fragmentedness replaces Baire class 1, is the working definition of tameness (cf. Theroem \ref{thm: old Prop 2}).

After that, we proceed to relate ``tameness of $(X,G)$'' / ``$\beta(X,G)<\infty$'' with one of the fundamental dividing lines in model theory: NIP. 
Research in this direction was already undertaken in model theory, e.g. \cite{Simon_2015}, \cite[Remark 1.6]{ArtemPierre}, in some sense also in \cite{New_bounded}.
Again, there were already existing results relating tameness of $(S_x(M),\aut(M))$ with NIP of $T$, where $M$ is a small model of $T$ (see Corollary 5.8 in \cite{KrRz}).
In fact, it was communicated to us after publishing the first version of this paper that 
the assumptions about being a small model in the results of Section 5 from \cite{KrRz} can be easily removed.
The results in Subsection \ref{subs: NIP characterization} which do not involve the $\beta$-rank
(so those concerning only NIP and tameness) 
either follow quickly from the results of Section 5 of \cite{KrRz} or from their proofs.
Nevertheless, 
our first main goal was to state explicite a relation between tameness of $(S_x(M),\aut(M))$
and NIP of $T$ (i.e. without any additional assumptions) and we provide 
all the proofs for a self-contained and more complete presentation.
Our second main goal here was to characterize NIP theories through our generalization of the $\beta$-rank. We achieve these both goals in Theorem \ref{thm: NIP, tame, rank} and Corollary \ref{cor: monster NIP, tame, rank}.
See also the variants for local space of types (Theorem \ref{thm: local NIP, tame, rank} and Corollary \ref{cor: new prop 6.6}).

In Section \ref{section: maps and bounds}, we discuss the relation between $\beta$-ranks of two dynamical systems $(X,G)$ and $(Y,H)$ such that there is a continuous surjection $\pi:X\to Y$ and corresponding maps between groups and Ellis semigroups. We naturally expect that the $\beta$-rank of an element from $E(Y,H)$ is bigger than or equal to the $\beta$-rank of its image in $E(X,G)$. However, this might be not true, except in the situation in which $\pi$ is an open map.
Finally, we sketch our further research plans related to application of the $\beta$-rank in model-theoretic context in Section \ref{section: postlude}.

We thank Tomasz Rzepecki for pointing out to us existing relations between stability of a theory $T$ and continuity of the elements of the Ellis semigroup, and relations between tameness of a dynamical system and NIP (see his doctoral dissertation \cite{rzepecki2018} and \cite{KrRz}). 
We also thank Krzysztof Krupiński for discussions that helped us understand the results of \cite{KrRz} better,
and Kyle Gannon, Anand Pillay and Ludomir Newelski for their valuable suggestions.
We express our gratitude towards the organizers of the GeSAMT for providing an excellent environment for communicating mathematical ideas.

\section{Definitions and first observations}\label{section: definitions and first observations}
Recall that $(X,G)$ is a dynamical system.
As $X$ is compact Hausdorff, we have that $X$ is uniformizable, and we denote by $\CU$ be the unique compatible uniformity, which consists of all the neighbourhoods of the diagonal in $X\times X$.

\begin{definition}
    Let $Y\subseteq X$ be closed, $f:Y\to X$, $y\in Y$ and $W\in\CU$. We say that $f$ $W$-oscillates at $y$ in $Y$ if
    for every $V\subseteq Y$, an open (in $Y$) neighbourhood of $y$, there exist $v_1, v_2\in V$ such that
    $(f(v_1),f(v_2))\not\in W$. For short, we write in this case $\osc(f,y)\geqslant W$.
\end{definition}

\begin{definition}
    Let $Y\subseteq X$ be closed, $f:Y\to X$ and $W\in\CU$. We define an $(f,W)$-derivative of the set $Y$ as
    $$(Y)'_{f,W}:=\{y\in Y\;|\;\osc(f,y)\geqslant W\}.$$
\end{definition}

\begin{remark}
In the notation of the definition above:
\begin{enumerate}
    \item 
    $(Y)'_{f,W}\subseteq Y$,
    \item 
    the $(f,W)$-derivative of $Y$ is again a closed subset of $X$.
    Too see that, we show that $Y\setminus (Y)'_{f,W}$ is open, so let $y\in Y\setminus (Y)'_{f,W}$.
    By the definition, there exists $V\subseteq Y$, an open neighbourhood of $y$ in $Y$, such that
    for every $v_1,v_2\in V$ we have $(f(v_1),f(v_2))\in W$. Clearly $y\in V\subseteq Y\setminus (Y)'_{f,W}$.
\end{enumerate}
\end{remark}

\begin{definition}
    Let $\alpha\in\Ord$, $Y\subseteq X$ be closed, $W\in\CU$ and let $f:Y\to X$.
    We define $(Y)^{\alpha}_{f,W}$ recursively as follows:
    \begin{itemize}
        \item $(Y)^{\alpha+1}_{f,W}:= \big( (Y)^{\alpha}_{f,W})'_{f|_{(Y)^{\alpha}_{f,W}},W}$,

        \item $(Y)^{\gamma}_{f,W}:=\bigcap\limits_{\alpha<\gamma} (Y)^{\alpha}_{f,W}$ for every limit ordinal $\gamma$.
    \end{itemize}
\end{definition}

\begin{definition}\label{definition: beta rank}
    \noindent\begin{enumerate}
        \item For a closed $Y\subseteq X$, $f:Y\to X$ and $W\in\CU$ we define the $\beta$-rank of the pair $(f,W)$ as $\beta(f,W):=\sup\{\alpha\in\Ord\;|\; (Y)^{\alpha}_{f,W}\neq\emptyset\}$. If $(Y)^\alpha_{f,W}\neq\emptyset$ for every $\alpha\in\Ord$ we say that $\beta(f,W)=\infty$.
    
        \item For $f:X\to X$ we define the $\beta$-rank of $f$ as $\beta(f):=\sup\{\beta(f,W)\;|\;W\in\CU\}$.

        \item Finally, we define the $\beta$-rank of the whole dynamical system as $\beta(X,G):=\sup\{\beta(\eta)\;|\;\eta\in E(X,G)\}$, where $E(X,G)$ denotes the Ellis semigroup of the action of $G$ on $X$.
    \end{enumerate}
\end{definition}

Let us comment here about the slight derivation from the original definition of the rank from \cite{GM}.
In the first item of the above definition, we define $\beta(f,W)$ as the supremum of $\alpha$'s such that $X^{\alpha}_{f.W}$ is non-empty. The original definition is a bit different and it takes $\beta(f,W)$ to be the least $\alpha$ such that $X^{\alpha}_{f,W}$ vanishes. Our reason here is to have the $\beta$-rank interact better with the model-theoretic Cantor-Bendixson rank on the space of types (e.g. p. 100 in \cite{tentzieg}): we set $X^0:=X$ and $X^{\alpha+1}$ to be $X^{\alpha}$ after removing isolated points,
and then the Cantor-Bendixson rank of $x\in X$ is the maximal $\alpha$ such that $x\in X^{\alpha}$ (see Lemma \ref{lemma: beta vs CB} for a precise statement). When $X$ is compact metric the unique compatible uniformity $\CU$ has as basis entourages of the form $U_\varepsilon=\{(x,y)\in X^2\mid d(x,y)<\varepsilon\}$, for $\varepsilon>0$. Using those entourages in the preceding definitions one recovers the notion of $\beta$-rank from \cite{GM}, except for the off-by-one difference we discussed in the previous paragraph.

\begin{remark}\label{remark: monotonicity}
    Let $f:X\to X$ and $W\in\CU$, and let $Y_1\subseteq Y_2$ be closed subsets of $X$.
    Then $(Y_1)'_{f,W}\subseteq (Y_2)'_{f,W}$. Assume that $y\in (Y_1)'_{f,W}\subseteq Y_1\subseteq Y_2$, and consider any open neighbourhood $V_2$ of $y$ in $Y_2$. We have that $V_1:=V_2\cap Y_1$ is an open neighbourhood of $y$ in $Y_1$ and by definition there exists $v_1,v_2\in V_1\subseteq V_2$ such that $(f(v_1),f(v_2))\not\in W$, hence $y\in (Y_2)'_{f,W}$.
    Using the proved inclusion, we conclude $\beta(f|_{Y_1},W)\leqslant \beta(f|_{Y_2},W)$ and $\beta(f|_{Y_1})\leqslant\beta(f_{Y_2})$. Moreover, the proved inclusion leads to $(Y_1)^{\alpha}_{f,W}\subseteq (Y_2)^{\alpha}_{f,W}$ for any $\alpha\in\Ord$.
\end{remark}

\begin{remark}\label{remark: subgroup}
    If $H\leqslant G$ then $\beta(X,H)\leqslant\beta(X,G)$.
\end{remark}

\begin{definition}
    Consider a closed subset $Y\subseteq X$, $f:Y\to X$, $W\in\CU$ and $x\in Y$.
    We define the following ranks:
    \begin{enumerate}
        \item $\bR_{f,W}(x)=\alpha$ if and only if $x\in (Y)^{\alpha}_{f,W}\setminus (Y)^{\alpha+1}_{f,W}$,
        \item $\bR_f(x)=\sup\limits_{U\in\CU}\bR_{f,U}(x)$,
        \item $\bR(x)=\sup\limits_{h\in E(X,G)}\bR_h(x)$.
    \end{enumerate}
\end{definition}

\begin{remark}\label{remark: continuity at point}
    A function $f:D\to X$, where $D\subseteq X$ is closed, is continuous at $d\in D$ if and only if
    $\bR_f(d)=0$. 
\end{remark}

\begin{proof}
First, let us note that $f$ is continuous at $d$ if and only if for every $W_2\in\CU$ there exists $W_1\in\CU$
such that $f[ W_1[d]\cap D]\subseteq W_2[f(d)]$.

Now, assume that $\bR_f(d)\neq 0$, so $1\leqslant\bR_f(d)=\sup\limits_{W\in\CU}\bR_{f,W}(d)$ and so there exists $W\in\CU$ such that $\bR_{f,W}(d)\geqslant 1$. Hence $d\in (D)'_{f,W}$, i.e. $d\in D$ and $\osc(f,d)\geqslant W$.
Take symmetric $W_2\in\CU$ such that $W_2\circ W_2\subseteq W$. Suppose that $f$ is continuous at $d$, then there exists $W_1\in\CU$ such that $f[W_1[d]\cap D]\subseteq W_2[f(d)]$.
Because $\osc(f,d)\geqslant W$ and $d\in W_1[d]\cap D$, which is open subset of $D$, we have $v_1,v_2\in W_1[d]\cap D$ such that $\big( f(v_1), f(v_2)\big)\not\in W$.
On the other hand, $f(v_1),f(v_2)\in W_2[f(d)]$ imply that $\big(f(v_1),f(v_2)\big)\in W_2\circ W_2\subseteq W$ - a contradiction.

If $f$ is not continuous at $d$, then there exist $W_2\in\CU$ such that for every $W_1\in\CU$ we have $f[W_1[d]\cap D]\not\subseteq W_2[f(d)]$. To see that $\osc(f,d)\geqslant W_2$, let $d\in V$ and let $V$ be an open subset of $D$. Then there exist $W_1\in\CU$ such that $W_1[d]\cap D\subseteq V$. We have $e\in W_1[d]\cap D$ with $f(e)\not\in W_2[f(d)]$, i.e $d,e\in V$ and $\big(f(d),f(e)\big)\not\in W_2.$
\end{proof}

\begin{lemma}\label{lemma: beta vs CB}
    Consider $f:X\to X$, $W\in\CU$ and $x\in X$.
    Then
    $$\bR_{f,W}(x)\leqslant \bR_f(x) \leqslant \bR(x) \leqslant \CB(x).$$
\end{lemma}

\begin{proof}
    First of all, let us notice that for any closed $D\subseteq X$, we have - by definitions - $(D)'_{f,W}\subseteq D'$. We claim that for every $\alpha\in\Ord$ it is $(X)^{\alpha}_{f,W}\subseteq X^\alpha$. The proof is by induction on $\alpha$. Case $\alpha=1$ follows from the first observation. Assume that $(X)^{\alpha}_{f,W}\subseteq X^{\alpha}$ and let us see what happens for $\alpha+1$.
    By Remark \ref{remark: monotonicity}, we get that $(X)^{\alpha+1}_{f,W}\subseteq (X^\alpha)'_{f,W}$.
    Then using the first observation with $D=X^{\alpha}$, we obtain $(X^{\alpha})'_{f,W}\subseteq X^{\alpha+1}$, and so $(X)^{\alpha+1}_{f,W}\subseteq X^{\alpha}$. The limit ordinal step is standard.

    Now, if $\bR_{f,W}(x)\geqslant\alpha$ then $x\in (X)^{\alpha}_{f,W}$, and so $x\in X^{\alpha}$ by the proved claim.
    The last thing implies that $\CB(x)\geqslant\alpha$. In other words, $\bR_{f,W}(x)\leqslant\CB(x)$, from which the lemma follows.
\end{proof}

\begin{cor}
    In the case of $X$ being the space of types $S_x(\FC)$ for a monster model $\FC\models T$ and a finite tuple of variables $x$, the Cantor-Bendixson rank of a point $p\in S_x(\FC)$ equals the Morley rank of $p$.
    For any topological group $G$ acting continuously on $S_x(\FC)$ via homeomorphisms (e.g. $\aut(\FC)$), and any 
    $f:S_x(\FC)\to S_x(\FC)$, $W\in\CU$ and $p\in S_x(\FC)$, we obtain
    $$\bR_{f,W}(p)\leqslant \bR_f(p) \leqslant \bR(p) \leqslant \MR(p).$$
\end{cor}

Let us mention here, that the $\beta$-rank behaves a bit different to the Morley rank. 
For example, if we consider some formula $\varphi$ and a monster model $\FC\models T$, we have
$$\MR(\varphi)=\max\{\MR(p)\;|\;\varphi\in p\in S(\FC)\}$$
(see Lemma 6.2.11 in \cite{tentzieg}). One could expect, that we can have a similar equality, i.e.
$\beta(f|_D,W)=\sup\{\bR_{f,W}(x)\;|\;x\in D\}$ for any closed $D\subseteq X$, $f:X\to X$ and any $W\in\CU$.
This is not the case - for example consider $d\in (X)'_{f,W}$ and take $D:=\{d\}$, then $\beta(f|_D,W)=0<\bR_{f,W}(d)$.
The reason behind it is that the $\beta$-rank takes into account the behaviour of $f$ on the open subsets of $D=\{d\}$, or the behaviour of $f$ on the open neighbourhoods of $d$ in $X$ in the second case. 
We accept this outcome of our definitions as in the general case of a closed $Y\subseteq X$, 
we are rather interested in the $Y$ as ``a dynamical sub-system'' and not just collection of points from $X$.
Anyway, we still obtain the following lemma.

\begin{lemma}
    Let $f:X\to X$, $W\in\CU$ and let $D\subseteq X$ be closed. Then
    $$\beta(f|_D,W)\leqslant\sup\{\bR_{f,W}(x)\;|\;x\in D\}.$$
\end{lemma}

\begin{proof}
    If $\beta(f|_D,W)\geqslant\alpha$ then $(D)^{\alpha}_{f,W}\neq\emptyset$. Choose $d\in (D)^{\alpha}_{f,W}\subseteq D$. By Remark \ref{remark: monotonicity}, $d\in (D)^{\alpha}_{f,W}\subseteq (X)^{\alpha}_{f,W}$ and so $\bR_{f,W}(d)\geqslant \alpha$.
\end{proof}

\begin{lemma}\label{lemma: image of derivative}
    Let $h:X\to X$ be a homeomorphism, $\alpha\in\Ord$, $f:X\to X$ and let $W\in\CU$, and let $Y\subseteq X$ be closed. Then
    $$h\big[ (Y)^{\alpha}_{f,W}\big] = (Y)^{\alpha}_{fh^{-1},W} = (Y)^{\alpha}_{hfh^{-1},h[W]}.$$
\end{lemma}

\begin{proof}
    The second equality is obvious directly from the definitions.
    For the proof of the first, we proceed with an induction. We start with the case of $\alpha=1$, so with showing that
    $$h\big[ (Y)'_{f,W}\big]= \big( h[Y] \big)'_{fh^{-1},W}.$$
    Let $x\in (Y)'_{f,W}$, then $h(y)\in h[Y]$. To show that $\osc(fh^{-1}|_{h[Y]},h(x))\geqslant W$, we take an open neighbourhood $U$ of $h(x)$ in $h[Y]$. Then $x\in h^{-1}[U]\subseteq Y$ is open, and - since $\osc(f|_Y,x)\geqslant W$ - there exist $v_1,v_2\in h^{-1}[U]$ such that $((f(v_1),f(v_2))\not\in W$. Then $h(v_1),h(v_2)\in U$ and $\big((fh^{-1})\big(h(v_1)\big),(fh^{-1})\big(h(v_2)\big)\big)\not\in W$. 
    Thus $h(x)\in \big( h[Y] \big)'_{fh^{-1},W}$. 
    As $h$, $Y$, $f$, $W$ are arbitrary, we replace them with $h^{-1}$, $h[Y]$, $fh^{-1}$ and $W$ respectively:
    $$h^{-1}\big[ (h[Y])'_{fh^{-1},W}\big]\subseteq \big( h^{-1}h[Y] \big)'_{fh^{-1}h,W}=\big(Y\big)'_{f,W}$$
    which, after applying $h$ to the both sides, gives us $(h[Y])'_{fh^{-1},W}\subseteq h[(Y)'_{f,W}]$ and so finally
    $h[(Y)'_{f,W}]=(h[Y])'_{fh^{-1},W}$. The rest of the induction is straightforward.
\end{proof}

\begin{cor}
    The rank $\bR(\cdot)$ is $G$-invariant, i.e. for every $h\in G$ and every $x\in X$ we have $\bR\big(h(x)\big)=\bR(x)$.
\end{cor}

\begin{proof}
    Let $h\in G$, $f:X\to X$, $x\in X$, $W\in\CU$ and let $\alpha\in\Ord$.
    If $\bR_{f,W}(x)\geqslant\alpha$ then $\bR_{fh^{-1},W}\big( h(x)\big)\geqslant\alpha$.
    To see this assume that $\bR_{f,W}(x)\geqslant\alpha$. 
    Then $x\in (X)^{\alpha}_{f,W}$ and $h(x)\in h\big[ (X)^{\alpha}_{f,W}\big] = (X)^{\alpha}_{fh^{-1},W}$ by Lemma \ref{lemma: image of derivative}, hence we get that $\bR_{fh^{-1},W}\big(h(x)\big)\geqslant\alpha$. By a symmetric argument, we get the opposite implication and therefore $\bR_{f,W}(x)=\bR_{fh^{-1},W}\big(h(x)\big)$.

    Note that $(E(X,G)h^{-1}=E(X,G)$ and so we can write
    \begin{IEEEeqnarray*}{rCl}
    \bR(x) & = & \sup\limits_{f\in E(X,G)}\sup\limits_{W\in\CU}\;\bR_{f,W}(x) = 
    \sup\limits_{f\in E(X,G)}\sup\limits_{W\in\CU}\;\bR_{fh^{-1},W}\big(h(x)\big) \\
    & = & \sup\limits_{f\in E(X,G)}\sup\limits_{W\in\CU}\;\bR_{f,W}\big(h(x)\big) \\
    & = & \bR\big(h(x)\big)
    \end{IEEEeqnarray*}
    which ends the proof.
    \end{proof}

The following is a corollary to Remark \ref{remark: continuity at point}, but we placed it here as it is more related to the last part of this section. The straightforward proof is left to the reader.

\begin{cor}\label{cor:beta0}
    For a dynamical system $(X,G)$, we have $\beta(X,G)=0$ if and only if every element of its Ellis semigroup is continuous.
\end{cor}

\begin{example}
    Consider the space of types $X:=S_x(\FC)$ in single variable $x$ for $\FC$ being a monster model of the theory ACF$_p$ ($p=0$ or $p$ being a prime number). Then, as a topological space $X$ has one non-isolated point (the type of a transcendental element) and the rest being isolated points (algebraic types). We claim that $\beta(X,G)=0$ for the natural action of $G:=\aut(\FC)$ on $X$.

    Let $\eta\in\EL(X,G)$, $W\in\CU$ be arbitrary. We are done if we show that $X'_{\eta,W}=\emptyset$.
    If $p$ is an isolated type, then obviously $\eta$ is continuous at $p$, so $\eta$ is not $W$-oscillating at $p$. 
    Let $p$ be the unique non-isolated type. We show that $\eta$ is continuous at $p$, which implies that is it is not $W$-oscillating at $p$.

    Let $\eta(p)\in [\psi]$ be a clopen neighbourhood of $\eta$. By compactness, $X\setminus[\psi]=\{q_1,\ldots,q_n\}$, where $q_1,\ldots,q_n$ are isolated.
    We have that $\eta$ may map isolated and non-isolated types to non-isolated ones, and $\eta$ maps the non-isolated type to itself. Moreover, if $q'$, $q^{\prime\prime}$ and $\eta(q')=\eta(q^{\prime\prime})$ are isolated, then $q'=q^{\prime\prime}$. Therefore $\eta^{-1}(\{q_1,\ldots,q_n\})$ is a finite set of isolated types and hence
    $$O:=[\psi]\;\cap\;(X\setminus \eta^{-1}(\{q_1,\ldots,q_n\})$$
    is open and we have $p\in O$ and $\eta(O)\subseteq [\psi]$.
\end{example}

\begin{prop} \label{prop: rank 0 iff stable}
    Let $\FC$ be a monster model of some $\CL$-theory $T$.
    The following are equivalent:
    \begin{enumerate}
        \item $T$ is stable,

        \item $\beta(S_x(\FC),\aut(\FC))=0$ for every finite tuple of variables $x$,

        \item $\beta(S_x(\FC),\aut(\FC))=0$ for every single variable $x$.
        
    \end{enumerate}
Moreover, by Corollary \ref{cor:beta0}, if for every single variable $x$, all elements of $E(S_x(\FC),\aut(\FC))$ are continuous, then the same holds for $E(S_x(\FC),\aut(\FC))$ where $x$ is any finite tuple of variables.
\end{prop}

\begin{proof}
    More or less, this is just binding together several already existing facts with our definitions, but let us explain how the proof goes exactly. Obviously there is nothing to prove in (2)$\Rightarrow$(3).

    One more observation before going to the actual proof.
    Fix some finite tuple of variables $x$ and set $(X,G):=(S_x(\FC),\aut(\FC))$.
    By Proposition II.2 in \cite{WAP} and by Corollary \ref{cor:beta0}, $\beta(X,G)=0$ if and only if $(X,G)$ is WAP (Weakly Almost Periodic), and the last means that each continuous complex-valued function $f$ on $X$ - denoted $f\in \Cc(X)$ - is WAP.

    Now, assume item (3). To show stability of $T$ it is enough to show that every formula $\varphi(x,y)$, where $|x|=1$, is stable. Take a formula $\varphi(x,y)$ with $|x|=1$.
    It follows that the characteristic function of a clopen set $[\varphi(x,a)]$, $\chi_{[\varphi(x,a)]}$ is WAP, which means that $\{\chi_{[\varphi(x,g^{-1}(a))]}\;|\;g\in G\}$ is weakly compact inside $\Cc_b(X)$ (bounded functions from $\Cc(x)$).
    By Grothendieck's theorem (e.g. Theorem 5.2 in \cite{WAPstable}), we have that for every sequence $(g_n)_{n<\omega}\subseteq G$  and every sequence $(p_m)_{m<\omega}\subseteq X$ it is
    $$\lim\limits_n \lim\limits_m \chi_{[\varphi(x,g^{-1}_n(a))]}(p_m)=\lim\limits_m \lim\limits_n \chi_{[\varphi(x,g^{-1}_n(a))]}(p_m).$$
    We can set $p_m:=\tp(b_m/\FC)$ for an arbitrary selected sequence of elements $(b_m)_{m<\omega}\subseteq\FC$, and so we get that for any sequence of elements $(a_n)_{n<\omega}$ sharing the same type over $\emptyset$,
    $$\lim\limits_n \lim\limits_m \varphi(b_m,a_n) = \lim\limits_m \lim\limits_n \varphi(b_m,a_n),$$
    which implies that there is no sequence $(b_m)_{m<\omega}\subseteq\FC$ and sequence $(a_n)_{n<\omega}$ such that
    $$\models\varphi(b_m,a_n)\quad\iff\quad m\leqslant n$$
    (the condition that $(a_n)_{n<\omega}$ share the same type over $\emptyset$ is removed by compactness).
    Thus $\varphi(x,y)$ is stable and so (1) follows.

    Now, assume (1), i.e. that $T$ is stable and fix some finite tuple of variables $x$.
    Then any formula $\varphi(x,y)$ is stable, and so
    for any sequences $(b_m)_{m<\omega}, (a_n)_{n<\omega}\subseteq\FC$ we have that 
    $$\lim\limits_n \lim\limits_m \varphi(b_m,a_n) = \lim\limits_m \lim\limits_n \varphi(b_m,a_n).$$
    Because of that, Theorem 5.2 in \cite{WAPstable} implies that $\chi_{[\varphi(x,a)]}$ is WAP for any $a\in\FC$.
    As $X$ is a zero-dimensional space, characteristic functions of clopen sets form a linearly dense subset of $\Cc(X)$.
    All WAP functions form a closed subalgebra in $\Cc(X)$ 
    (see the paragraph after Fact 2.1 in \cite{Iba16} for the case of real valued functions $C(X)$, then consult \cite{Ben13}, Fact 2), 
    hence we see that every function from $\Cc(X)$ is WAP, so $(S_x(\FC), \aut(\FC))$ is WAP and we obtain item (2).
\end{proof}

\noindent
Let us note here, that above Proposition \ref{prop: rank 0 iff stable} can be also proved using Proposition 2.2 from \cite{Anand_Gro}.


\section{o-minimal theories}\label{section: o-minimal}

The goal of this section is to prove that if $T$ is an o-minimal theory (in some language $\mathcal L$ including an order $<$), then for all models $M$ of $T$, we have 
$$\beta(S_x(M),\aut(M))\leqslant 1,$$
where $x$ is a single variable. Together with Proposition \ref{prop: rank 0 iff stable}, we can conclude $\beta(S_x(\FC),\aut(\FC))=1$ whenever $\FC$ is a monster model of an o-minimal theory.
Until the end of this section, fix an o-minimal theory $T$, a model $M$ of $T$ and let $X:=S_x(M)$ (where $|x|=1$), $G:=\aut(M)$.
We have a natural order on $X$ given by 
\begin{IEEEeqnarray*}{rCl}
p(x)<^* q(x) &\iff & \text{there exists }m\in M\text{ such that }\\
& & (``x<m"\in p \text{ and }``m\leqslant x"\in q )\\
&\text{or}&  (``x\leqslant m"\in p\text{ and }``m<x"\in q).
\end{IEEEeqnarray*}
The order $<^*$ is a strict total order on $X$. The proof is standard, the only thing which needs some argument is to show that for every distinct $p,q\in X$ we have $p<^* q$ or $q<^* p$. Assume that $p\neq q$ but not $p<^*q$ nor $q<^*p$.
Then for every $m\in M$, all the following points hold at the same time:
\begin{itemize}
    \item $``x\geqslant m"\in p$ or $``m>x"\in q$,
    \item $``x>m"\in p$ or $``m\geqslant x"\in q$,
    \item $``x\geqslant m"\in q$ or $``m>x"\in p$,
    \item $``x>m"\in q$ or $``m\geqslant x"\in p$.
\end{itemize}
Using the above, we conclude for example that: $``x\geqslant m"\in p$ if and only if $``x\geqslant m"\in q$;
$``x>m"\in p$ if and only if $``x>m"\in q$. As $p\neq q$ and $X$ is Hausdorff, there exists $\psi(x)\in\mathcal{L}$
such that $p\in[\psi]$ and $q\not\in[\psi]$. As $T$ is o-minimal, $\psi$ can be written as a finite union of points in $M$
and intervals with endpoints in $M$. But then $p\in[\psi]$ implies that $q\in[\psi]$, a contradiction.

\begin{remark} \label{remark: basic stuff on the o-minimal case}
 For every $f\in E(X,G)$ and $p,q\in X$, if $p<^* q$ then $f(p)\leqslant^* f(q)$ (i.e $f(p)<^*f(q)$ or $f(p)=f(q)$).
\end{remark}

\begin{proof}
    Let $a\models p$ and $b\models q$ (in some monster model $\FC$). As $p<^* q$, we have that $a<b$.
    If not $f(p)\leqslant^* f(q)$, then $f(q)<^*f(p)$ and so there exist $n\in M$ such that, without loss of generality,
    $f(q)\in[x\leqslant n]$ and $f(p)\in[n<x]$. If $(\eta_i)_{i\in I}$ is a net in $G$ pointwise converging to $f$,
    then there exists $i_0\in I$ such that for every $i\geqslant i_0$, we have $\eta_i(q)\in[x\leqslant n]$ and $\eta_i(p)\in[n<x]$. Take $i\geqslant i_0$, then $q\in[x\leqslant\eta^{-1}_i(n)]$ and $p\in[\eta^{-1}_i(n)<x]$.
    So $b\leqslant\eta^{-1}_i(n)<a$ which contradicts $a<b$.
\end{proof}

 Our goal is to prove that for any $f\in E(X,G)$ and any $W\in\CU$ (the unique compatible uniformity on $X$), $(X)^2_{f,W}=\emptyset$. We will show that $(X)^1_{f,W}$ is finite, from which the conclusion follows immediately. The arguments are very close to those of \cite[Section 4]{codenotti2023}, slightly modified to work when $X$ is not necessarily metrizable.

For any $W\in\CU$ fix $U(W)\in\CU$ symmetric with $U(W)\circ U(W)\subseteq W$ and a symmetric set $V\in\CU$ with $V\circ V\subseteq U(W)$. For every $p\in X$ consider a clopen interval $p\in I_p=[a_p\leqslant x\leqslant b_p]\subseteq V[p]$ ($a_p,b_p\in M$ or $I_p$ is a half-interval, i.e. an interval without one endpoint and with one endpoint in $M$).
Then $(I_p)_{p\in X}$ forms an open cover of $X$ and so there exists $n(W)<\omega$ and intervals $I^W_1,\ldots,I^W_{n(W)}$
chosen among $(I_p)_{p\in X}$ which also form an open cover of $X$.

\begin{definition}
    Let $p\in X$, $f\in E(X,G)$, $W\in\CU$ and $i\leqslant n(W)$ be such that $f(p)\in I^W_i$. We say that $f$ is $U(W)$-oscillating in the right direction at $p$ if for every open neighbourhood $V$ of $p$, there exists $p^V\in V$ such that
    \begin{itemize}
        \item $(f(p),f(p^V))\not\in U(W)$,
        \item $r<^*f(p^V)$ for every $r\in I^W_i$.
    \end{itemize}
    If the second condition holds with $f(p^V)<^*r$ for every $r\in I^W_i$ instead we say that $f$ is $U(W)$-oscillating in the left direction at $p$. 
    Note that it is possible for $f$ to oscillate in both the left and right direction at some $p\in X$.
\end{definition}

\begin{lemma}\label{lemma: f must oscillate in some direction}
    Let $f\in E(X,G)$, $W\in\CU$, $i\leqslant n(W)$ and $p\in X$ be such that $f(p)\in I^W_i$ and $p\in (X)^1_{f,W}$. Then $f$ must be $U(W)$-oscillating in at least one direction at $p$. 
\end{lemma}

\begin{proof}
    Since $p\in (X)^1_{f,W}$, for every open neighbourhood $V$ of $p$, we can find $r^V,q^V\in V$ with $(f(r^V),f(q^V))\not\in W$. Since $U(W)$ satisfies $U(W)\circ U(W)\subseteq W$, we cannot have both $(f(r^V),f(p))\in U(W)$ and $(f(p),f(q^V))\in U(W)$. If the first fails, let $p^V:=r^V$, otherwise let $p^V:=q^V$. In any case we have $(f(p),f(p^V))\not\in U(W)$. Since $f(p)\in I^W_i\subseteq V[p']$ and $V\circ V\subseteq U(W)$ (for some choice of $p'\in X$ and $V\in\CU$, as explained earlier) - $f(p^V)\in I^W_i$ would lead to $(f(p),f(p^V))\in V\circ V\subseteq U(W)$, so to a contradiction. Thus $f(p^V)\not\in I^W_i$.
    Recall that $I^W_i=[a\leqslant x\leqslant b]$ for some $a,b\in M$. Thus $f(p^V)\not\in I^W_i$ implies that $``a>x"\in f(p^V)$ or $``x>b"\in f(p^V)$.
    
    The last thing means that either $r<^*f(p^V)$ for every $r\in I^W_i$, or that $f(p^V)<^*r$ for all $r\in I^W_i$.
    As we iterate over all open neighbourhoods of $p$ one of those inequalities holds infinitely many times by the pigeonhole principle, showing that $f$ is $U(W)$-oscillating in some direction at $p$.
\end{proof}

\begin{lemma}\label{lemma: f cannot oscillate in the same direction at different points}
    If $p_1<^*p_2\in (X)^1_{f,W}$ are such that, for some $i\leqslant n(W)$, $f(p_1)\in I^W_i$ and $f(p_2)\in I^W_i$,
    then $f$ cannot be $U(W)$-oscillating in the same direction at $p_1$ and $p_2$.
\end{lemma}

\begin{proof}
    Suppose for a contradiction that $f$ is $U(W)$-oscillating in the same direction at $p_1$ and $p_2$, which we may assume to be the right direction without loss of generality. 
    As $p_1<^* p_2$ there is $m\in M$ such that $``x<m"\in p_1$ and $``m\leqslant x"\in p_2$, or $``x\leqslant m"\in p_1$ and $``m<x"\in p_2$. If $``x<m"\in p_1$, then we set $V:=[x<m]$, and let $V:=[x\leqslant m]$ in the other case.
    We see that $p_1\in V$ and $p_2\not\in V$.

    Let $p_1^V\in V$ witness that $f$ is $U(W)$-oscillating in the right direction at $p_1$, i.e. $(f(p_1),f(p^V_1))\not\in U(W)$ and $r<^* f(p^V_1)$ for any $r\in I^W_i$. As $p_1^V\in V$, we have that $p_1^V<^*p_2$.
    By Remark \ref{remark: basic stuff on the o-minimal case}, we get $f(p_1^V)\leqslant^* f(p_2)$ which contradicts $f(p_2)\in I^W_i$.
%
\end{proof}

\begin{theorem}\label{thm: o-min in x}
    Let $f\in E(X,G)$. Then $(X)^1_{f,W}$ is finite for all $W\in\CU$. In particular $\beta(f)\leqslant 1$ and $\beta(X,G)\leqslant 1$.
\end{theorem}

\begin{proof}
    Let $W\in\CU$. By Lemma \ref{lemma: f must oscillate in some direction} if $p\in (X)^1_{f,W}$ and $f(p)\in I^W_i$, then $f$ must $U(W)$-oscillate in some direction at $p$. By Lemma \ref{lemma: f cannot oscillate in the same direction at different points} (and again Lemma \ref{lemma: f must oscillate in some direction}), there can be at most one another point $q\in(X)^1_{f,W}$ with $f(q)\in I^W_i$. 
    As there are finitely many $I^W_i$, this shows that $(X)^1_{f,W}$ is finite, hence discrete. In particular $(X)^2_{f,W}=\emptyset$. Since $W\in\CU$ was arbitrary we have $\beta(f)\leqslant 1$, and since $f\in E(X,G)$ was arbitrary we have $\beta(X,G)\leqslant 1$.
\end{proof}

It can be checked that with the same strategy as in the proof above, one can obtain the following variant of Theorem \ref{thm: o-min in x}.

\begin{theorem}[Order topology]
    Let $(X,<)$ be a dense strict total order such that the order topology induced on $X$ is Hausdorff and compact.
    Moreover, let $G$ be a topological group continuously acting on $X$ by order preserving bijections.
    Then $(X)^1_{f,W}$ is finite for all $W\in\CU$. In particular $\beta(f)\leqslant 1$ and $\beta(X,G)\leqslant 1$.
\end{theorem}

In particular we obtain the result stated at the beginning of this section

\begin{theorem}\label{theorem: rank of o-minimal theories}
    Let $T$ be an o-minimal theory, $M$ a model of $T$ and $x$ a single variable. Then $\beta(S_x(M),\aut(M))\leqslant 1$.
\end{theorem}

\begin{question}
    The above shows that $\beta(S_x(\FC),\aut(\FC))=1$ when $\FC$ is a monster model of some o-minimal theory and $|x|=1$. What happens when $x$ is a bigger tuple?
\end{question}

\subsection{A counterexample}
In view of Theorem \ref{theorem: rank of o-minimal theories} it is natural to ask for an example of a theory $T$ which is not o-minimal, but with the property that $\beta(S_x(M),\aut(M))= 1$ whenever $M\vDash T$ and $x$ is a single variable. The goal of this section is to show that a dense cyclic order provides such an example. The arguments will be less detailed, being similar to those employed in the o-minimal case.

\begin{definition}
    Let $L=\{C\}$ be a language consisting of a single ternary relation and let $M$ be the structure with underlying set $\mathbb Q$ and $C$ interpreted as the cyclic order $$C(a,b,c)\iff (a<b<c)\lor(b<c<a)\lor(c<a<b).$$
    Let $T=\Th(\mathbb Q,C)$, which is a NIP theory, but it is not o-minimal, nor stable.
\end{definition}

For the rest of this section fix $M\vDash T$, $X=S_x(M)$, where $x$ is a single variable and $G=\aut(M)$. Our goal is to show that $\beta(X,G)\leqslant 1$ (as $T$ is not stable, we already have $\beta(X,G)>0$ by Proposition \ref{prop: rank 0 iff stable}, provided $M$ is a monster model).

Similarly to the o-minimal case we can extend $C$ to a relation $C^\ast$ on $X$, where for $p_1,p_2,p_3\in X$ we have \begin{align*}
  C^\ast(p_1,p_2,p_3)\iff (\exists a_1,a_2,a_3\in M)&\big(M\models C(a_1,a_2,a_3)\land ``C(a_1,x,a_2)"\in p_1\land\\ &``C(a_2,x,a_3)"\in p_2\land ``C(a_3,x,a_1)"\in p_3\big).
\end{align*}
As in the proof of Remark \ref{remark: basic stuff on the o-minimal case}, one can check that, for every $f\in E(X,G)$, if $C^\ast(p_1,p_2,p_3)$ holds, then either $C^\ast(f(p_1),f(p_2),f(p_3))$ holds, or $f$ is not injective on $\{p_1,p_2,p_3\}$.

Once again, instead of showing that $(X)^2_{f,W}=\emptyset$ for every $f\in E(X,G)$ and every $W\in\CU$ (the unique compatible uniformity of $X$), we will show that $(X)^1_{f,W}$ is finite, hence discrete, from which $\beta(f,W)\leqslant 1$ follows immediately.

For any $W\in\CU$ fix $U(W)\in\CU$ symmetric with $U(W)\circ U(W)\subseteq W$ and a symmetric set $V\in\CU$ with $V\circ V\subseteq U(W)$. For every $p\in X$ consider a clopen set $p\in I_p$ determined by the formula $C(a,x,b)$ for some $a,b\in M$, so that $I_p\subseteq V[p]$. Since $(I_p)_{p\in X}$ is an open cover of $X$, there exists $n(W)<\omega$ and sets $I^W_1,\ldots,I^W_{n(W)}$ chosen among $(I_p)_{p\in X}$ which still cover $X$. 

Arguing as in the beginning of the proof of Lemma \ref{lemma: f must oscillate in some direction}, we see that if $p\in (X)^1_{f,W}$, for some $f\in E(X,G)$ and some $W\in\CU$, then for every open neighbourhood $U$ of $p$, we can find $p^U\in U$ such that $(f(p),f(p^U))\not\in U(W)$.
\begin{lemma}\label{lemma: finitely many pi in cyclic order}
    Let $p_1,p_2,p_3\in X$ be distinct and $f\in E(X,G)$ be such that $f(p_j)\in I^W_i$ for $j\in 1,2,3$ and a fixed $i\leqslant n(W)$. Then at least one of the $p_j$ is not in $(X)^1_{f,W}$.
\end{lemma}
\begin{proof}
    Assume without loss of generality that $C^\ast(p_1,p_2,p_3)$ holds. Suppose for a contradiction that $p_j\in (X)^1_{f,W}$ for every $j\in 1,2,3$ and fix an open neighbourhood $U$ of $p_2$ such that $U\cap\{p_1,p_2,p_3\}=p_2$. By assumption we can find $p_2'\in U$ such that $(f(p_2),f(p_2'))\not\in V$. Note that we cannot have $f(p_1)=f(p_3)$, because otherwise, as $f$ preserves the nonstrict version of $C^\ast$, $f$ would have to be constant on the set determined by the formula $"C(p_1,x,p_3)"$, which contains $p_2$, contradicting our assumption that $p_2\in (X)^1_{f,W}$. Moreover $f(p_2')\neq f(p_1)$ and $f(p_2')\neq f(p_3)$, since $f(p_1),f(p_2),f(p_3)\in I^W_i$, but $I^W_i\circ I^W_i\subseteq V$ and $(f(p_2),f(p_2'))\not\in V$.  By construction we have $C^\ast(p_1,p_2',p_3)$, so we must have $C^\ast(f(p_1),f(p_2'),f(p_3))$, since $f(p_1)\neq f(p_2')\neq f(p_3)$. However since $f(p_1),f(p_2)\in I^W_i$, while $f(p_2')\not\in I^W_i$, we have $C^\ast(f(p_1),f(p_2),f(p_2'))$ instead, a contradiction.
\end{proof}

\begin{theorem} Let $f\in E(X,G)$. Then $(X)^1_{f,W}$ is finite for all $W\in\CU$. In particular $\beta(f)\leqslant 1$ and $\beta(X,G)\leqslant 1$.
\end{theorem}
\begin{proof} By Lemma \ref{lemma: finitely many pi in cyclic order} if $p\in(X)^1_{f,W}$ and $p\in I^W_i$, there can be at most another point $q\in(X)^1_{f,W}$ with $f(q)\in I^W_i$. Since $X$ is covered by finitely many sets of the form $I^W_i$, this shows that $(X)^1_{f,W}$ is finite, hence discrete. In particular $(X)^2_{f,W}=\emptyset$ and $\beta(f,W)\leqslant 1$. Since $W\in\CU$ was arbitrary we have $\beta(f)\leqslant 1$, and since $f\in E(X,G)$ was arbitrary we have $\beta(X,G)\leqslant 1$.    
\end{proof}

\begin{question}
    Can we characterize an interesting dividing line in the stability hierarchy by the condition $\beta(S_x(\FC),\aut(\FC))=1$ for every single variable $x$?
\end{question}

\section{Tame systems and NIP theories}\label{section: NIP theories} 
By Proposition \ref{prop: rank 0 iff stable}, we know that stability of $T$ is equivalent to the fact that every element of $E(X,G)$ is a Baire class $0$ function. In this section, we are interested in the situation when $T$ is NIP, which is somehow related to the fact that every element of $E(X,G)$ is a Baire class $1$ function. Let us study this problem and relate it to the $\beta$-rank.

\subsection{Generalizing the ingredients}
First, we provide a theorem which can be proved using already existing tools. Then we will generalize these tools. Let us have a look at the theorem and its proof.

\begin{theorem}\label{thm: NIP for small}
    Let $T$ be a theory in a countable language. Let $M$ be a countable model of $T$. Then $T$ is NIP if and only if for all finite tuples $x$ and all $f\in E(S_x(M),\aut(M))$ we have $\beta(f)<\omega_1$
\end{theorem}
\begin{proof}
  By \cite[Corollary 5.8]{KrRz}, $T$ is NIP if and only if $(S_x(M),\aut(M))$ is tame for all finite tuples $x$. Since $M$ and $T$ are countable, $S_x(M)$ is metrizable, so by \cite[Theorem 6.3]{GlasnerMegrelishviliUspenskij}, the tameness of $(S_x(M),\aut(M))$ is equivalent to the fact that every $f\in E(S_x(M),\aut(M))$ is a Baire class 1 function $X\to X$, which in turn is equivalent to $\beta(f)<\omega_1$ by \cite[essentially Proposition 2]{KechrisLouveau}.
\end{proof}

In the above proof the $\beta$-rank of $f$ is related to $f$ being a Baire class $1$ function via Proposition 2 from \cite{KechrisLouveau}, which is stated in the case of metrizable $X$. Theorem \ref{thm: old Prop 2} is our generalization of Proposition 2 from \cite{KechrisLouveau} to the case of a general (compact Hausdorff) space $X$. First of all, we obtain different bounds, but this is natural. Second, instead of being a Baire class $1$ function, we consider being a \emph{fragmented} function. In the case of $X$ being a Polish space, $f:X\to X$ is fragmented if and only if $f$
is barely continuous, if and only if $f$ is Baire class $1$ (see Lemma 2.2 in \cite{GM23}).

Therefore, the equivalence between first and second item in Theorem \ref{thm: old Prop 2} coincides with Proposition 2 in \cite{KechrisLouveau} for $X$ being a Polish space.

\begin{definition}[Definition 2.1 in \protect{\cite{GM23}}]
    A function $f:X\to X$ is fragmented if 
    for every non-empty closed $D\subseteq X$
    for every $W\in\CU$ there exists an open $V\subseteq X$ such that
    $D\cap V\neq\emptyset$ and $f[D\cap V \times D\cap V]\subseteq W$.
\end{definition}

\begin{lemma}\label{lemma: fragmented vs derivative}
    Let $f:X\to X$. The following are equivalent.
    \begin{enumerate}
        \item $f$ is fragmented.
        
        \item For every non-empty closed set $D\subseteq X$, $f$ has arbitrarily small oscillation on $D$.
        That is, for every non-empty closed set $D\subseteq X$, for every $W\in\CU$ there exists $d\in D$ such that
        $\osc(f|_D,d)<W$, i.e. there exists an open neighbourhood $V$ of $d$ in $D$ such that $f[V\times V]\subseteq W$.

        \item For every non-empty closed $D\subseteq X$, for every $W\in\CU$ we have $(D)'_{f,W}\subsetneq D$.
    \end{enumerate} 
  Moreover, if $f$ is barely continuous then $f$ is fragmented.  
\end{lemma}

\begin{proof}
    The proof is straightforward from the definitions.
\end{proof}

Recall that the weight of a topological space $X$, denoted by $w(X)$ is defined to be the smallest possible cardinality of a basis of $X$. In particular if $X$ is a compact metric space, $w(X)=\aleph_0$.

\begin{theorem}\label{thm: old Prop 2}
    Let $f:X\to X$. The following are equivalent.
    \begin{enumerate}
        \item $f$ is fragmented,
        \item $\beta(f)<w(X)^+$,
        \item $\beta(f)<\infty$.
    \end{enumerate}
\end{theorem}

\begin{proof}
  (1)$\Rightarrow$(2): Assume that $f$ is fragmented and let $W\in\CU$.
    By Lemma \ref{lemma: fragmented vs derivative}, if $D\subseteq X$ is non-empty and closed, then there exists $d\in D$ such that $d\in D\setminus (D)'_{f|_D,W}$.

    We claim that for every ordinals $\alpha<\beta$ we have $(X)^{\beta}_{f,W}\subsetneq (X)^{\alpha}_{f,W}$. 
    To see this take $\alpha\in\Ord$. Because $(X)^{\alpha}_{f,W}$ is a closed subset of $X$, we can apply the first paragraph to $D=(X)^{\alpha}_{f,W}$ to get $(X)^{\alpha+1}_{f,W}=(D)'_{f|_D,W}\subsetneq D=(X)^{\alpha}_{f,W}$. 
    Thus on every step when the rank $\beta(f,W)$ increases, we throw out at least one point of $X$.
    In other words, the chain $\big((X)^{\alpha}_{f,W}\big)$ stabilizes at $\emptyset$ after $\lambda$ steps for some $\lambda<|X|^+$.
    
    Let $\kappa=w(X)$ and let $(U_\alpha)_{\alpha<\kappa}$ be an enumeration of a basis of $X$. We prove that if $(F_{\beta})_{\beta<\lambda}$ is a strictly decreasing chain of closed sets in $X$, then $\lambda<\kappa^+$. 
    To every $F_\beta$ we associate $C_\beta\subseteq\kappa$ defined by 
    $$C_\beta=\{\alpha\in\kappa\mid U_\alpha\cap F_\beta\neq\emptyset\}.$$ 
    Since $X\setminus F_\beta=\bigcup_{\alpha\not\in C_\beta} U_\alpha$, we have that the map $F_\beta\mapsto C_\beta$ is injective. Moreover if $F_\beta\subsetneq F_\gamma$, then $C_\beta\subsetneq C_\gamma$. In other words $(C_\beta)_{\beta<\lambda}$ is a strictly decreasing well ordered chain of subsets of $\kappa$, hence $\lambda<\kappa^+$. 
    To prove $(2)$ it only remains to observe that, $\big((X)^{\alpha}_{f,W}\big)_{\alpha<\lambda}$ is strictly decreasing (because $f$ is fragmented) chain of closed sets in $X$, and so 
    $\beta(f,W)<w(X)^+$.

    Observe that for any closed set $D\subseteq X$ and $W_1,W_2\in\CU$, if $W_1\subseteq W_2$ then $(D)'_{f,W_1}\supseteq (D)'_{f,W_2}$. From that we obtain that $(X)^{\alpha}_{f,W_1}\supseteq (X)^{\alpha}_{f,W_2}$ whenever $W_1\subseteq W_2$. Finally, $W_1\subseteq W_2$ implies that $\beta(f,W_1)\geqslant \beta(f,W_2)$.

    Now, recall that $W\in\CU$, so 
    $$W=\bigcup\limits_{i\in I, j\in J}U_i\times U_j,$$
    for some collection of indices $I$ and $J$. As the diagonal $\Delta\subseteq X\times X$ is a compact set, there exists finite $I_0\subseteq I$ and finite $J_0\subseteq J$ such that
    $$\Delta\subseteq W_0:=\bigcup\limits_{i\in I_0, j\in J_0}U_i\times U_j\subseteq W.$$
    Then $\beta(f,W)\leqslant\beta(f,W_0)$.
    To count $\beta(f)=\sup\limits_{W\in\CU}\beta(f,W)$, we can take the supremum over sets of the form as $W_0$ above, instead over the whole $\CU$. Then the supremum is taken over a set of size $w(X)$ and, as $w(X)^+$ is regular, it must be that $\beta(f)<w(X)^+$.

    (2)$\Rightarrow$(3) is obvious. Let us move to (3)$\Rightarrow$(1). 
    Assume that $f$ is not fragmented. By Lemma \ref{lemma: fragmented vs derivative}, there exists a non-empty closed $D\subseteq X$ and $W\in\CU$ such that $(D)'_{f,W}=D$.
    This is sufficient to show that $\beta(f)=\beta(f,W)=\infty$, as $\emptyset\neq D=(D)^{\alpha}_{f,W}\subseteq (X)^{\alpha}_{f,W}$ for any $\alpha\in\Ord$.
\end{proof}

The second ingredient of the proof of Theorem \ref{thm: NIP for small}, i.e. Theorem 6.3 from \cite{GlasnerMegrelishviliUspenskij}, is removed by replacing ``Baire class $1$'' with ``fragmented''.
More precisely, we borrow the following definition from \cite{GM23}.

\begin{definition}[Definition 3.1 in \cite{GM23}]
    The dynamical system $(X,G)$ is tame if every $f\in E(X,G)$ is fragmented.
\end{definition}

\begin{cor}\label{cor: tame vs defined rank}
    The following are equivalent:
    \begin{enumerate}
        \item The dynamical system $(X,G)$ is tame.
        \item $\beta(X,G)<\infty$.
        \item For every $f\in E(X,G)$, $\beta(f)<w(X)^+$
    \end{enumerate} 
\end{cor}
\begin{proof}
    This is clear by Theorem \ref{thm: old Prop 2}.
\end{proof}

Recall that we work all the time with $(X,G)$ being a dynamical system, where $X$ is Hausdorff compact.
Let us note here, that in the case of $X$ being additionally metric, Theorem 3.1 from \cite{GM23} combined with Theorem 6.3 from \cite{GlasnerMegrelishviliUspenskij} say that $(X,G)$ is tame if and only if every $f\in E(X,G)$ is Baire class $1$, if and only if $(X,G)$ is tame in the sense of the old definition of tameness (appearing in the proof of Theorem \ref{thm: NIP for small}). Now, we are going to investigate a generalization of the last ingredient of the proof of Theorem \ref{thm: NIP for small}, i.e. Corollary 5.8 from \cite{KrRz}.

\begin{definition}
    We say that a sequence $(f_n)_{n<\omega}\subseteq C(X)$ is independent if there exist real numbers $r_1<r_2$ such that for every finite and disjoint $P,M\subseteq\mathbb{N}$ we have
    $$\emptyset\neq\bigcap\limits_{n\in P}f_n^{-1}[(-\infty,r_1)]\;\cap\;\bigcap\limits_{n\in M}f_n^{-1}[(r_2,\infty)].$$
\end{definition}

Note that by Theorem 2.4 from \cite{GM23}, a bounded subset $F\subseteq C(X)$ does not contain an independent sequence if and only if it does not contain an $\ell^1$-sequence.

\begin{definition}
    Let $f\in C(X)$. We say that $f$ is tame in $(X,G)$ if one of the following equivalent conditions holds.
    \begin{enumerate}
        \item $\{f\circ g\;|\;g\in G\}$ does not contain an independent sequence.
        \item $\{f\circ g\;|\;g\in G\}$ does not contain an $\ell^1$-sequence.
    \end{enumerate}
\end{definition}

\begin{remark}\label{remark: tame system vs tame function}
    The dynamical system $(X,G)$ is tame if and only if every $f\in C(X)$ is tame in $(X,G)$ (see Definition 3.1 in \cite{GM23}).
\end{remark}

\subsection{Characterization of NIP}\label{subs: NIP characterization}
In this subsection fix a language $\CL$ and a complete $\CL$-theory $T$ with a monster model $\FC$,
which strongly $\kappa$-homogeneous and $\kappa$-saturated. A model $M\models T$ is \emph{small} if $|M|<\kappa$.

\begin{definition}
    Let $\varphi(x,y)$ be an $\CL$-formula.
    \begin{enumerate}
        \item We say that $\varphi(x,y)$ is small-tame if for every small $M\models T$ and every $b\in M^y$
        the characteristic function $\chi_{[\varphi(x,b)]}$ is tame in $(S_x(M),\aut(M))$ (cf. Definition 4.20 and Definition 2.69 in \cite{rzepecki2018}, or Definition 5.3 in \cite{KrRz}).

        \item We say that $\varphi(x,y)$ is tame in $M\models T$ if for every $b\in M^y$
        the characteristic function $\chi_{[\varphi(x,b)]}$ is tame in $(S_x(M),\aut(M))$.
        
        \item We say that $\varphi(x,y)$ is tame if it is tame in $M$ for every $M\models T$.
    \end{enumerate}
\end{definition}

\begin{lemma}\label{lemma: tame vs NIP}
    Let $\varphi(x,y)$ be an $\CL$-formula and let $N\models T$ be strongly $\aleph_0$-homogeneous and $\aleph_1$-saturated.
    The following are equivalent.
    \begin{enumerate}
        \item $\varphi(x,y)$ is NIP,
        \item $\varphi(x,y)$ is tame,
        \item $\varphi(x,y)$ is small-tame,
        \item $\varphi(x,y)$ is tame in $N$.
    \end{enumerate}
\end{lemma}

\begin{proof}
    (2)$\Rightarrow$(3) and (2)$\Rightarrow$(4) follow by definition. Equivalence between (1) and (3) was proved in Lemma 5.4 in \cite{KrRz}. We only show (1)$\Rightarrow$(2) and (4)$\Rightarrow$(2).

    Assume that $\varphi(x,y)$ is not tame. There exists $M\models T$ and $b\in M^y$ such that $\chi_{[\varphi(x,b)]}$ is not tame, i.e. there is a sequence $(g_n)_{n<\omega}\subseteq\aut(M)$ such that
    $$\chi_{[\varphi(x,b)]}\circ g_n=\chi_{[\varphi(x,g_n(b))]}$$
    forms an independent sequence.

    Recall that $\varphi(x,y)$ has IP if there exists a sequence $(b_n)_{n<\omega}\subseteq\FC$ such that
    $\big(\varphi(\FC,b_n)\big)_{n<\omega}$ is an independent family, i.e. for every $I\subseteq\mathbb{N}$
    $$\emptyset\neq \bigcap\limits_{n\in I}\varphi(\FC,b_n)\;\cap\;\bigcap\limits_{n\in\mathbb{N}\setminus I}\neg\varphi(\FC,b_n).$$
    We set $b_n:=g_n(b)$. Since $(\chi_{[\varphi(x,b_n)]})_{n<\omega}$ is an independent sequence, there exist real numbers $r_1<r_2$ such that for every disjoint finite $P_+,P_-\subseteq\mathbb{N}$
    $$\emptyset\neq\bigcap\limits_{n\in P_-}\chi^{-1}_{[\varphi(x,b_n)]}[(-\infty,r_1)]\;\cap\;\bigcap\limits_{n\in P_+}\chi^{-1}_{[\varphi(x,b_n)]}[(r_2,\infty)].$$
    As $P_+$ and $P_-$ are arbitrary it follows that we can take $r_1=\frac{1}{3}$ and $r_2=\frac{2}{3}$.
    To see that $\big(\varphi(\FC,b_n)\big)_{n<\omega}$ is an independent family, consider arbitrary $I\subseteq\mathbb{N}$ and set $J:=\mathbb{N}\setminus I$ and 
    $$\pi(x):=\{\varphi(x,b_i),\,\neg\varphi(x,b_j)\;|\;i\in I,\,j\in J\}.$$
    We are done if we show that $\pi(x)$ is consistent, which we achieve by compactness. Let $P_+\subseteq I$ and $P_-\subseteq J$ be finite and let $p\in S_x(M)$ be such that
    $$p\in \bigcap\limits_{n\in P_-}\chi^{-1}_{[\varphi(x,b_n)]}[(-\infty,\frac{1}{3})]\;\cap\;\bigcap\limits_{n\in P_+}\chi^{-1}_{[\varphi(x,b_n)]}[(\frac{2}{3},\infty)].$$
    Then $\{\neg\varphi(x,b_j),\,\varphi(x,b_i)\;|\;j\in P_-,\,i\in P_+\}\subseteq p$.

    To show (4)$\Rightarrow$(2), we will shot that not (2) implies not (4). If there is $M\models T$ and $b\in M^y$ such that $\chi_{[\varphi(x,b)]}$ is not tame in $(S_x(M),\aut(M))$, we can find $(f_n)_{n<\omega}\subseteq\aut(M)$
    such that for every finite $P,M\subseteq\Nn$ with $P\cap M=\emptyset$, the set
    $$\{\neg\varphi(x,f_i(b)),\,\varphi(x,f_j(b))\;|\;i\in P,\in j\in M\}$$
    is consistent. Let $b_n:=f_n(b)$, $\bar{b}:=(b_n)_{n<\omega}$. Let $\bar{c}\models\tp(\bar{b})$ be a tuple contained in $N$.   
    Then for every $n<\omega$, there exists $g_n\in\aut(N)$ such that $c_n=g_n(c_0)$ (since $b_n$ and $b_0$ have the same type over $\emptyset$ and $N$ is strongly $\aleph_0$-homogeneous). Thus for every finite $P,M\subseteq\Nn$ with $P\cap M=\emptyset$, the set
     $$\{\neg\varphi(x,g_i(c_0)),\,\varphi(x,g_j(c_0))\;|\;i\in P,\in j\in M\}$$
    is consistent. This contradicts tameness of $\chi_{[\varphi(x,c_0)]}$ in $(S_x(N),\aut(N))$ and $\varphi(x,y)$ is not tame in $N$.
\end{proof}

\noindent
In the above proof, one could adapt the proof of the ``moreover'' part from Corollary 5.9 in \cite{KrRz} to show (4)$\Rightarrow$(1) instead of showing (4)$\Rightarrow$(2), but it would not shorten the proof significantly.
On the other hand, we can use (1)$\iff$(3) (i.e. Lemma 5.4 from \cite{KrRz}) to show (1)$\Rightarrow$(2).
More precisely, if $\varphi(x,y)$ is not tame in some $M$, this model $M$ can be embedded in an even bigger monster model $\FC'$, such that $M$ is small in $\FC'$. Then (1)$\iff$(3) considered for $\FC'$ implies that $\varphi(x,y)$ is IP. We are leaving the above proof of Lemma \ref{lemma: tame vs NIP} as it is more straightforward and shows the interactions between the main definitions.

In the next theorem, the equivalence between (1) and (2) appears in some form in Corollary 5.8 in \cite{KrRz} - it will be the equivalence (1)$\iff$(3) from Corollary 5.8 in \cite{KrRz}, but after noticing that the assumption of ``smallnes'' there can be removed.

\begin{theorem}\label{thm: NIP, tame, rank}
    The following are equivalent.
    \begin{enumerate}
        \item $T$ is NIP.
        
        \item For every $M\models T$ and every finite tuple of variables $x$, the dynamical system $(S_x(M),\aut(M))$ is tame.
        
        \item For every $M\models T$ and every finite tuple of variables $x$, 
        $$\beta\big(S_x(M),\aut(M) \big)<\infty.$$

        \item For every $M\models T$ and every finite tuple of variables $x$, if $f\in E(S_x(M),\aut(M))$ then $\beta(f)<w(S_x(M))^+ \leqslant (|M|+|T|)^+$.

        \item For every $M\models T$ and every single variable $x$, the dynamical system \\ $(S_x(M),\aut(M))$ is tame.
        
        \item For every $M\models T$ and every single variable $x$, 
        $$\beta\big(S_x(M),\aut(M) \big)<\infty.$$

        \item For every $M\models T$ and every single variable $x$, if $f\in E(S_x(M),\aut(M))$ then $\beta(f)<w(S_x(M))^+ \leqslant (|M|+|T|)^+$.        
    \end{enumerate}
\end{theorem}

\begin{proof}
    The equivalences between (2), (3) and (4) follow by Corollary \ref{cor: tame vs defined rank}.
    Similarly we obtain equivalences between (5), (6) and (7). Then, we can note that (4)$\Rightarrow$(5).
    
    For (5)$\Rightarrow$(1), we argue as follows. By Lemma \ref{lemma: tame vs NIP}, an $\CL$-formula $\varphi(x,y)$ is NIP
    if and only if for any model $M\models T$ and any $b\in M^y$,
    the characteristic function $\chi_{[\varphi(x,y)]}$ is tame in $(S_x(M),\aut(M))$. 
    To show that $T$ is NIP, it is enough to show that every formula $\varphi(x,y)$ with $|x|=1$ is NIP (by Proposition 2.11 and Lemma 2.5 from \cite{Guide_NIP}).
    The last thing follows by (5) and Remark \ref{remark: tame system vs tame function}.

    Let us move to (1)$\Rightarrow$(2). First, select a model $M\models T$, we want to show that $(S_x(M),\aut(M))$ is tame, i.e.
    that every $f\in C(S_x(M))$ is tame in $(S_x(M),\aut(M))$. Note that as $S_x(M)$ is a zero-dimensional space, the characteristic functions of clopen sets form a linearly dense subset of $C(S_x(M))$.
    By Lemma \ref{lemma: tame vs NIP}, the characteristic functions of clopen sets are tame.
    Now, as tame functions form a closed subalgebra of $C(S_x(M))$ (see the third paragraph after Corollary 3.2 in \cite{Iba16}) we obtain that all functions in $C(S_x(M))$ are tame.
\end{proof}

In the above theorem we obtained a characterization of NIP in the terms of the $\beta$-rank and of tameness of the dynamical system $(S_x(M),\aut(M))$.
The authors of \cite{ArtemPierre} noted in their Remark 1.6 that the NIP assumption implies tameness of the dynamical system $(S_G(M),G(M))$ (where $G$ is a definable group), but it is not equivalent to it. Using our Corollary \ref{cor: tame vs defined rank}, one can conclude that $\beta(S_G(M),G(M))<\infty$.

\begin{remark}
Another tame dynamical system involving a definable group $G$ in the NIP context, appears in Remark 6.11 in \cite{Artem_Kyle}. Let $T$ be a NIP theory, in a countable language $\CL$, expanding a group. 
Consider monster model $\CG\models T$ and a countable $G\preceq\CG$. By $S_x(\CG,G)$ we denote the space of types from $S_x(\CG)$ which are finitely satisfiable in $G$, similarly $\mathfrak{M}_x(\CG,G)$ is the space of Keisler measures finitely satisfiable in $G$. First, we note that the dynamical system $(S_x(\CG,G),G)$ is tame (this was communicated to us by Kyle Gannon), thus $\beta(S_x(\CG,G),G)<\infty$. Then, we use Remark 6.11 from \cite{Artem_Kyle} to conclude that also $\beta(\mathfrak{M}_x(\CG,G),G)<\infty$.
\end{remark}

\begin{question}
    Suppose that $T$ is a countable NIP theory. 
    If $M$ and $M'$ are two countable models of $T$, can we have $$\beta(S_x(M),\aut(M))\neq\beta(S_x(M'),\aut(M'))\text{?}$$
\end{question}

\begin{theorem}\label{thm: one vs many}
    Let $x$ be a finite tuple of variables and let $N\models T$ be strongly $\aleph_0$-homogeneous and $\aleph_1$-saturated.
    The following are equivalent.
\begin{enumerate}
    \item $(S_x(N),\aut(N))$ is tame,

    \item for every $M\models T$, $(S_x(M),\aut(M))$ is tame.
\end{enumerate}
\end{theorem}

\begin{proof}
    Assume that there exists $M\models T$ such that $(S_x(M),\aut(M))$ is not tame.
    By Remark \ref{remark: tame system vs tame function}, there exists $f\in C(S_x(M))$ which is not tame in $(S_x(M),\aut(M))$. If all the characteristic functions of all the clopen sets in $S_x(M)$ would be tame, then every element of $C(S_x(M))$ would be a tame function. Thus we can assume that $f=\chi_{[\varphi(x,b)]}$ for some $\CL$-formula $\varphi(x,y)$ and $b\in M^y$. Hence, $\varphi(x,y)$ is not tame, and - by Lemma \ref{lemma: tame vs NIP} - $\varphi(x,y)$ is not tame in $N$.
\end{proof}
\noindent
Note that the above theorem is similar to the "moreover part" in Corollary 5.9 from \cite{KrRz} (and can be derived from it combined with Lemma \ref{lemma: tame vs NIP}).
The next corollary will be important later, for the development of local ranks as it assures us that we can work in one previously selected monster model instead of considering all models of a given theory.

\begin{cor}\label{cor: monster NIP, tame, rank}
    Let $\FC\models T$ be a monster model. The following are equivalent.
    \begin{enumerate}
        \item $T$ is NIP.

        \item For every single variable $x$, the dynamical system $(S_x(\FC),\aut(\FC))$ is tame.

        \item For every single variable $x$, $\beta(S_x(\FC),\aut(\FC))<\infty$.

        \item For every single variable $x$, if $f\in E(S_x(\FC),\aut(\FC))$ then $\beta(f)<w(S_x(\FC))^+ $.
    \end{enumerate}
\end{cor}

\begin{proof}
    (2), (3) and (4) are equivalent by Corollary \ref{cor: tame vs defined rank}. (1)$\Rightarrow$(2) follows by Theorem \ref{thm: NIP, tame, rank}. (2)$\Rightarrow$(1) follows by Theorem \ref{thm: one vs many} and Theorem \ref{thm: NIP, tame, rank}.
\end{proof}

\subsection{Local variants}\label{subs: local}
Here, we restate the main results from the previous subsection, but for the space of $\varphi$-types, or more generally for the space of $\Delta$-types (see below). We keep proofs short, as they are similar to the previous ones, but for clarity we still include them.

Let $T$ be a complete $\CL$-theory with a monster model $\FC$, i.e. $\FC\models T$ and $\FC$ is strongly $\kappa$-homogeneous and $\kappa$-saturated for some big cardinal $\kappa$. Fix a finite tuple of variables $x$
and a set of $\CL$-formulas
$$\Delta(x):=\{\varphi_i(x,y_i)\;|\;i\in I\}.$$
We allow $\Delta(x)$ to be all $\CL$-formulas in variables $x$ - then we recover the general setting from the previous part.
If $A\subseteq\FC$, then by $S_{\Delta(x)}(A)$ we denote the space of $\Delta$-types, i.e. the set of maximal consistent collections of Boolean combinations of instances of formulas from $\Delta(x)$. 
If $M\preceq\FC$, then we can consider the dynamical system $(S_{\Delta(x)}(M),\aut(M))$, where $\aut(M)$ acts on $S_{\Delta(x)}(M)$ in the natural way. To prove the following remark it is enough to expand the definitions and use (model-theoretic) compactness, so we skip the proof.

\begin{remark}\label{rem: local tame}
    Let $\varphi(x,y)$ be a Boolean combination of formulas from $\Delta(x)$, $M\models T$ and $b\in M^y$.
    Then $\chi_{[\varphi(x,b)]}\in C(S_x(M))$ is tame in $(S_x(M),\aut(M))$ if and only if
    $\chi_{[\varphi(x,b)]}\in C(S_{\Delta(x)}(M))$ is tame in $(S_{\Delta(x)}(M),\aut(M))$.
\end{remark}

\noindent
The following theorem is a local variant of Theorem \ref{thm: NIP, tame, rank} (we fixed the tuple of variables, thus we have less points).

\begin{theorem}\label{thm: local NIP, tame, rank}
    The following are equivalent.
\begin{enumerate}
    \item Every $\varphi\in\Delta$ is NIP.

    \item $(S_{\Delta(x)}(M),\aut(M))$ is tame for every $M\models T$.

    \item $\beta(S_{\Delta(x)}(M),\aut(M))<\infty$ for every $M\models T$.

    \item $\beta(f)<w(S_{\Delta(x)}(M))^+$ for every $M\models T$ and every $f\in E(S_{\Delta(x)}(M),\aut(M))$.
\end{enumerate}
\end{theorem}

\begin{proof}
    As previously, equivalences between (2), (3) and (4) are established via Corollary \ref{cor: tame vs defined rank}.
    To prove (2), we need to show that every $f\in C(S_{\Delta(x)}(M))$ is tame in $(S_{\Delta(x)}(M),\aut(M))$, where $M\models T$. It is enough to check the last condition for the characteristic functions of clopen sets in $S_{\Delta(x)}(M)$, thus let $\varphi(x,y)$ be a Boolean combination of formulas from $\Delta(x)$.
    Point (1) implies that $\varphi$ is NIP.
    By Lemma \ref{lemma: tame vs NIP}, $\varphi$ is tame in $M$,
    i.e. for every $b\in M^y$, $\chi_{[\varphi(x,b)]}$ is tame in $(S_x(M),\aut(M))$. 
    Then, Remark \ref{rem: local tame}, gives us that every such $\chi_{[\varphi(x,b)]}$ is tame in $(S_{\Delta(x)}(M),\aut(M))$. The proof of the implication from (2) to (1) uses the same ingredients.
\end{proof}

\begin{lemma}\label{lemma: local one vs many}
    Let $N\models T$ be strongly $\aleph_0$-homogeneous and $\aleph_1$-saturated. The following are equivalent.
    \begin{enumerate}
        \item $(S_{\Delta(x)}(N),\aut(N))$ is tame,
        \item for very $M\models T$, $(S_{\Delta(x)}(M),\aut(M))$ is tame.
    \end{enumerate}
\end{lemma}

\begin{proof}
    If there exists $M\models T$, a Boolean combination $\varphi(x,y)$ of formulas from $\Delta(x)$ and
    $b\in M^y$ such that $\chi_{[\varphi(x,b)]}$ is not tame in $(S_{\Delta(x)}(M),\aut(M))$, then, by Remark \ref{rem: local tame}, we see that $\chi_{[\varphi(x,b)]}$ is not tame in $(S_{x}(M),\aut(M))$.
    As $\varphi(x,y)$ is not tame, Lemma \ref{lemma: tame vs NIP} imples that $\varphi(x,y)$ is not tame in $N$.
    This, by Remark \ref{rem: local tame}, results in untameness of $(S_{\Delta(x)}(N),\aut(N))$.
\end{proof}

By combining Theorem \ref{thm: local NIP, tame, rank} with Lemma \ref{lemma: local one vs many}, we obtain the following corollary, which is analogous to Proposition 6.6 from \cite{casasimpl}. However, to be closer to the statement of Proposition 6.6 from \cite{casasimpl}, a notion of a new local rank is needed - see Section \ref{section: postlude}.

\begin{cor}\label{cor: new prop 6.6}
    The following are equivalent.
    \begin{enumerate}
        \item Every $\varphi\in\Delta$ is NIP.

        \item $(S_{\Delta(x)}(\FC),\aut(\FC))$ is tame.

        \item $\beta(S_{\Delta(x)}(\FC),\aut(\FC))<\infty$.

        \item $\beta(f)<w(S_{\Delta(x)}(\FC))^+$ for every $f\in E(S_{\Delta(x)}(\FC),\aut(\FC))$.
    \end{enumerate}
\end{cor}

\subsection{Examples of NIP theories with higher rank}
    In this section we show examples of NIP theories whose models have $\beta$-rank bigger than one. 
    After introducing the notion of $\beta$-rank in \cite{GM}, Glasner and Megrelishvili pointed out that all the examples of tame dynamical systems known at the time had $\beta$-rank at most $1$ (or rank $2$ in their notation, see the comment after Definition \ref{definition: beta rank}) and asked \cite[Question 11.8]{GM} for examples with higher $\beta$-rank. While examples of arbitrary $\beta$-rank have already been constructed in \cite[Section 5]{codenotti2023}, our Lemma \ref{lemma: rank is at least n} give a new family of examples.

\begin{definition}
    Let $\mathcal{L}_n:=\{<_1,\ldots,<_n\}$ and 
    $T_n$ be the $\CL_n$-theory $\Th(\mathbb Q^n,<_1,\ldots,<_n)$, where for $x=(x_1,\ldots,x_n)$, $y=(y_1,\ldots,y_n)$, we interpret 
    $$x <_i y\iff x_i<y_i.$$ 
\end{definition}

Fix some $n<\omega$. 
The theory $T_n$ is NIP but not o-minimal. To see that $T_n$ is NIP but not o-minimal, we note that dp-rank of $T_n$ is equal to $n$ and every o-minimal theory is dp-minimal, i.e. has dp-rank bounded by $1$. Then as dp-rank of $T_n$ is finite it must be NIP. Using standard techniques, it might be shown that $T_n$ has quantifier elimination.

Choose $M\models T_n$ and let us work for a moment in $S_x(M)$ where $x$ is a single variable.
Let $(\xi_1,\ldots,\xi_n)\in \{-,+\}^{\times n}$. We define
$$p^0_{(\xi_1,\ldots,\xi_n)}(x):=\bigcup\limits_{\substack{i\leqslant n\\ \text{ such that } \\  \xi_i=+}}\{m<_i x\;|\;m\in M\}\;\cup\;
\bigcup\limits_{\substack{i\leqslant n\\ \text{ such that } \\  \xi_i=-}}\{x<_i m\;|\;m\in M\}.$$
By quantifier elimination, for every choice of $(\xi_1,\ldots,\xi_n)$, $p^0_{(\xi_1,\ldots,\xi_n)}(x)$ determines a complete type, which we denote by $p_{(\xi_1,\ldots,\xi_n)}(x)\in S_x(M)$.

\begin{lemma}\label{lemma: rank is at least n}
    Let $M:=(\mathbb Q^n,<_1,\ldots,<_n)$, $G:=\aut(M)$ and $X=S_x(M)$, where $|x|=1$. We have $$\beta(X,G)\geqslant n.$$
\end{lemma}

\begin{proof}
    The proof will be in two steps. First we will describe a function $f\colon X\to X$ and prove that $\beta(f)=n$. Afterward we will prove that $f\in E(X,G)$, concluding the proof.

    In order to describe $f$ we need to introduce some terminology. 
    We say that $p(x)\in S_x(M)$ is at $+\infty$ (respectively at $-\infty$) according to $i\leqslant n$ if for every $m\in M$, the formula ``$m<_ix$'' belongs to $p(x)$ (respectively the formula ``$x<_i m$'' belongs to $p(x)$).
    
    For $p(x)\in X$ we define 
    $$f(p(x))=p_{(\xi_1,\ldots,\xi_n)}(x),$$
    where $\xi_i=-$ if and only if $p(x)$ is at $-\infty$ according to $i$, and $\xi_i=+$ otherwise. 
    Now let $W\in\CU$ (the unique compatible uniformity on $X$) be such that 
    $$(p_{\xi_1,\ldots,\xi_n)},p_{(\zeta_1,\ldots,\zeta_n)})\not\in W$$
    for any distinct $(\xi_1,\ldots,\xi_n),(\zeta_1,\ldots,\zeta_n)\in\{-,+\}^{\times n}$.
    We will show that $\beta(f,W)=n$. 
    We begin by computing $(X)^1_{f,W}$:

    If $p(x)\in X$ is such that for every $i\leqslant n$,
    $p(x)$ is not at $-\infty$ according to $i$, then $f$ is locally constant at $p(x)$ (with value the type $p_{(+,\ldots,+)}(x)$), 
    so there is no $W$-oscillation of $f$ at $p(x)$.

    If instead $p(x)\in X$ is such that for some $i\leqslant n$, $p(x)$ is at $-\infty$ according to $i$, then every neighbourhood of $p(x)$ must contain a type $q(x)$ such that for every $j\leqslant n$, $q(x)$ is not at $-\infty$ according to $j$. 
    Then $f(p(x))=p_{(\xi_1,\ldots,\xi_n)}(x)$ and $f(q(x))=p_{(\zeta_1,\ldots,\zeta_n)}(x)$ are such that
    $\xi_i=-$ and $\zeta_i=+$, so that, by our choice of $W$, $\big(f\big(p(x)\big),f\big(q(x)\big)\big)\not\in W$.

    The above shows that $(X)^1_{f,W}$ consists exactly of the types that are at $-\infty$ according to at least one $i\leqslant n$. 
    Iterating the same argument we obtain that $(X)^j_{f,W}$, for $j\leqslant n$, consists of all the types that are at $-\infty$ according to at least $j$ distinct $i$'s. In particular $(X)^n_{f,W}$ contains only the type which is at $-\infty$ according to every $i\leqslant n$, while $(X)^{n+1}_{f,W}=\emptyset$, showing that $\beta(f,W)=n$. Since the same argument works verbatim for any $W'\in\CU$ with $W'\subseteq W$, this shows that $\beta(f)=n$.

    It only remains to show that $f\in E(X,G)$, but it is easy to check that $f$ is the pointwise limit of $(g_i)_{i\in\mathbb N}$, where $g_i((x_1,\ldots,x_n))=(x_1+i,\ldots,x_n+i)$.
\end{proof}

\begin{question}
    How is the $\beta$-rank related to the dp-rank in the class of NIP theories?
\end{question}

\section{Maps and bounds}\label{section: maps and bounds}

Below, we state a slight generalization of Theorem 7 \cite[p. 54]{Auslander} - a standard fact, but we could not find it in the literature. The proof is straightforward, so we skip it.

\begin{lemma}\label{lemma: two dynamical systems}
Let $(X,G)$ and $(Y,H)$ be dynamical systems, let $\rho: G\to H$ be a continuous epimorphism and $\pi:X\to Y$ a continuous surjection such that $\rho(g)\pi(x)=\pi(gx)$ for every $g\in G$ and $x\in X$.
Then there exists a unique map $\theta:E(X,G)\to E(Y,H)$ such that
$$\theta(f)\big(\pi(x)\big)=\pi \big( f(x)\big)$$
for every $f\in E(X,G)$ and every $x\in X$.
Moreover, such map $\theta$ is a surjective continuous semigroup homomorphism.
\end{lemma}

To the end of this section, let us work with the assumptions of the above lemma, i.e.:
let $(X,G)$ and $(Y,H)$ be dynamical systems, let $\rho: G\to H$ be a continuous epimorphism and $\pi:X\to Y$ a continuous surjection such that $\rho(g)\pi(x)=\pi(gx)$ for every $g\in G$ and $x\in X$.
Moreover, let $\theta:E(X,G)\to E(Y,H)$ be the map from the thesis of the lemma.

\begin{lemma}\label{lemma: derived sets vs preaimage}
Consider any closed $D\subseteq Y$, $f\in E(X,G)$ and $U\in\CU_Y$ (the uniformity of $Y$), and $\alpha\in \Ord$.
We have
\begin{enumerate}
    \item $$\big(\pi^{-1}[D]\big)^{\alpha}_{f,\pi^{-1}[U]}\subseteq\pi^{-1}\big[ (D)^{\alpha}_{\theta(f),U}\big],$$
    \item if moreover $\pi$ is an open map then 
    $$\big(\pi^{-1}[D]\big)^{\alpha}_{f,\pi^{-1}[U]} = \pi^{-1}\big[ (D)^{\alpha}_{\theta(f),U}\big].$$
\end{enumerate}
\end{lemma}

\begin{proof}
    We prove the first point recursively starting with the case of $\alpha=1$. For that, let us recall that $x\in \pi^{-1}\big[(D)'_{\theta(f),U}\big]$ if and only if:
    \begin{itemize}
        \item $\pi(x)\in D$ and
        \item for every open neighbourhood $W$ of $\pi(x)$ in $D$ there exist $y_1,y_2\in W$ such that 
        $\big( \theta(f)(y_1),\theta(f)(y_2)\big)\not\in U$.
    \end{itemize}
    On the other hand, after using properties of $\theta$, we can write that 
    $x\in \big( \pi^{-1}[D]\big)'_{f,\pi^{-1}[U]}$ if and only if
    \begin{itemize}
        \item $\pi(x)\in D$ and
        \item for every open neighbourhood $B$ of $x$ in $\pi^{-1}[D]$ there exist $x_1,x_2\in B$ such that
        $\Big( \theta(f)\big(\pi(x_1)\big), \theta(f)\big(\pi(x_2)\big) \Big)\not\in U$.
    \end{itemize}
Now, if $x\in \big( \pi^{-1}[D]\big)'_{f,\pi^{-1}[U]}$ and $W$ is an open neighbourhood of $\pi(x)$ in $D$,
then $\pi^{-1}[W]$ is an open neighbourhood of $x$ in $\pi^{-1}[D]$ and so we obtain $x_1,x_2\in\pi^{-1}[W]$
such that for $y_1:=\pi(x_1), y_2:=\pi(x_2)\in W$ we have $\big( \theta(f)(y_1),\theta(f)(y_2)\big)\not\in U$.
Therefore 
$$\big(\pi^{-1}[D]\big)'_{f,\pi^{-1}[U]}\subseteq\pi^{-1}\big[ (D)'_{\theta(f),U}\big].$$

Assume the first point and let us show that its counterpart holds for $\alpha+1$.
By Remark \ref{remark: monotonicity}, the induction hypothesis and the case of $\alpha=1$ for arbitrary closed set $D$, we get the following:
\begin{IEEEeqnarray*}{rCl}
\big( \pi^{-1}[D]\big)^{\alpha+1}_{f,\pi^{-1}[U]} & = & \Big( \big(\pi^{-1}[D]\big)^{\alpha}_{f,\pi^{-1}[U]}\Big)'_{f,\pi^{-1}[U]}  \subseteq \Big( \pi^{-1}\big[(D)^{\alpha}_{\theta(f),U} \big]\Big)'_{f,\pi^{-1}[U]} \subseteq \\
& \subseteq & \pi^{-1}\Big[\big( (D)^{\alpha}_{\theta(f),U}\big)'_{\theta(f),U} \Big] = \pi^{-1}\big[(D)^{\alpha+1}_{\theta(f),U} \big]
\end{IEEEeqnarray*}
The step with a limit ordinal is standard and omitted.

To show the second point of the lemma, we repeat the above strategy of the proof after achieving 
$$\big(\pi^{-1}[D]\big)'_{f,\pi^{-1}[U]} = \pi^{-1}\big[ (D)'_{\theta(f),U}\big].$$
To obtain the above, we need to start with $x\in \pi^{-1}\big[ (D)'_{\theta(f),U}\big]$ and an open neighbourhood $B$ of $x$ in $\pi^{-1}[D]$. As $\pi$ is an open map, we have that $\pi[B]$ is an open neighbourhood of $\pi(x)$ in $D$.
Thus there exist $y_1,y_2\in \pi[B]$ such that $\big( \theta(f)(y_1),\theta(f)(y_2)\big)\not\in U$.
As $y_1,y_2\in \pi[B]$, there exist $x_1,x_2\in B$ such that $y_1=\pi(x_1)$ and $y_2=\pi(x_2)$, and so 
$\Big( \theta(f)\big(\pi(x_1)\big), \theta(f)\big(\pi(x_2)\big) \Big)\not\in U$. This ends the proof.
\end{proof}

\begin{lemma}\label{lemma: rank in two systems}
Consider any $f\in E(X,G)$ and $U\in\CU_Y$ (the uniformity of $Y$). We have
\begin{enumerate}
    \item $\beta(f,\pi^{-1}[U])\leqslant \beta(\theta(f),U)$,
    \item if moreover $\pi$ is an open map, then
    $\beta(f,\pi^{-1}[U])=\beta(\theta(f),U)$ and so $\beta(f)\geqslant\beta(\theta(f))$.
\end{enumerate}
\end{lemma}

\begin{proof}
    For the proof, simply use Lemma \ref{lemma: derived sets vs preaimage} and the definitions.
\end{proof}

The fact that, in Lemma \ref{lemma: rank in two systems}, we have $\beta(f,\pi^{-1}[U])\leqslant \beta(\theta(f),U)$
is a bit counter-intuitive as the potentially smaller system could have a bigger rank. What will happen in the extreme case of $Y=\{y\}$, $H=\{ 1\}$ and $(X,G)$ being a system with $\beta(X,G)>0$?
Of course, $E(Y,H)=\{\id_Y\}$, hence $\beta(\theta(f),U)=0$ for any $f\in E(X,G)$ and any $U\in \CU_Y$, and so 
$\beta(f,\pi^{-1}[U])=0$ by the above lemma. However, there is no choice in picking up $U\in\CU_Y$, i.e. $U=\{y\}\times\{y\}$ and so $\pi^{-1}[U]=X\times X$, and $\beta(f,X\times X)=0$ for obvious reasons.
This example shows that the relation between ranks of $(X,G)$ and $(Y,H)$ is more subtle and depends on the relation between the uniformities on $X$ and $Y$. Note that in our case, the continuous map $\pi$ is automatically uniformly continuous, so restricting our attention to the uniformly continuous maps does not help.

\begin{example}
    Assume that $M\models T$ and $x'$ is a subtuple of a finite tuple of variables $x$.
    Then the restriction map $\pi:S_{x}(M)\to S_{x'}(M)$ is an open, surjective map and we can use Lemma \ref{lemma: rank in two systems}, to conclude that
    $$\beta(S_x(M),\aut(M))\geqslant\beta(S_{x'}(M),\aut(M)).$$
\end{example}

\section{Postlude}\label{section: postlude}
Consider a monster model $\FC\models T$ of some theory $T$.
Let $H:=\{\sigma\in\aut(\FC)\;|\;\sigma(M)=M\}$, where $M\preceq\FC$.
Then, by Remark \ref{remark: subgroup}, we have 
$$\beta(S_x(\FC),\aut(\FC))\geqslant\beta(S_x(\FC),H).$$
The group $H$ acts on $S_x(\FC)$ in such a way that the restriction maps $\pi: S_x(\FC)\to S_x(M)$ and $\rho:H\to \aut(M)$ commute with corresponding group actions, so we are in the situation of Lemma \ref{lemma: two dynamical systems}. If $\pi$ would be an open map, we could conclude that 
$$\beta(S_x(\FC),\aut(\FC))\geqslant \beta(S_x(\FC),H)\geqslant\beta(S_x(M),\aut(M)).$$
Unfortunately the map $\pi$ is rarely an open map. Moreover, if $q\in S_x(\FC)$ then it does not always follow that $\bR(q)\geqslant\bR(\pi(q))$. Preserving the Cantor-Bendixson rank is related to non-forking extensions of types, thus we ask the following natural question.

\begin{question}
    What can we say about a type extensions $p|_M\subseteq p\in S_x(\FC)$ if it preserves the $\beta$-rank?
\end{question}

\noindent
The question above is a bit vague, namely what does it mean to preserve the $\beta$-rank?
Let us mimic the definitions of (local) ranks (and multiplicities)
derived from the Cantor-Bendixson rank on $S(\FC)$ which were used to characterize type extensions in the case of stable theories (eg. as in Chapter 6 of \cite{casasimpl}).
Similarly as in Subsection \ref{subs: local}, we fix some finite tuple of variables $x$ and a set of formulas
$\Delta(x)=\{\varphi_i(x,y_i)\;|\;i\in I\}$. Moreover, let $\pi(x)$ be a set of $\CL(A)$-formulas, where $A\subseteq\FC$.
As the map
$$\res_{\Delta}:S_x(\FC)\to S_{\Delta(x)}(\FC)$$
is continuous and a closed map, the following set is closed
$$X_{\pi,\Delta}:=\res_{\Delta}\Big[\bigcap\limits_{\psi(x)\in\pi(x)}[\psi(x)]\Big]=\{p(x)\in S_{\Delta(x)}(\FC)\;|\;p(x)\,\cup\,\pi(x)\text{ is consistent}\}.$$

\begin{definition}
Let $f\in E(S_{\Delta(x)}(\FC),\aut(\FC))$, $W\in\CU$. We define
    $$\beta^{\Delta}_{f,W}(\pi):=\beta(f|_{X_{\pi,\Delta}},W).$$

\end{definition}

Now, we can reformulate the previous question to make it precise:

\begin{question}
    What can we say about a type extensions $p|_M\subseteq p\in S_x(\FC)$ if
    for every finite set of formulas $\Delta(x)$, every $f\in E(S_{\Delta(x)}(\FC),\aut(\FC))$ and each $W\in \CU$,
    $$\beta^{\Delta}_{f,W}(p|_M)=\bR_{f|_{X_{\pi,\Delta}},W}(p|_{\Delta})?$$
\end{question}

As it was earlier promised, we rephrase Corollary \ref{cor: new prop 6.6} so its statement 
looks more similar to the one of Proposition 6.6 in \cite{casasimpl}.

\begin{cor}\label{cor: new prop 6.6 vol.2}
    The following are equivalent.
    \begin{enumerate}
        \item Every $\varphi\in\Delta$ is NIP.

        \item $\beta^{\Delta}_{f,W}(x=x)<w\big(S_{\Delta(x)}(\FC)\big)^+$ for every $f\in E(S_{\Delta(x)}(\FC),\aut(\FC))$ and every $W\in\CU$.

        \item $\beta^{\Delta}_{f,W}(x=x)<\infty$ for every $f\in E(S_{\Delta(x)}(\FC),\aut(\FC))$ and every $W\in\CU$.
    \end{enumerate}
\end{cor}

Finally, let us mention that in \cite{KechrisLouveau}, there were defined 3 ranks: $\alpha$, $\beta$ and $\gamma$. It would be worth to understand the interactions of the remaining ranks $\alpha$ and $\gamma$ with model theory.
For example, note that in the definition of set $B^{\S}(X)$ in \cite{Simon_2015} (Definition 2.1 and Theorem 2.3 there)
there is a similar idea as in the definition of the separation rank $\alpha$ from \cite{KechrisLouveau}.

\printbibliography
\end{document}